\documentclass[12pt, reqno]{amsart}

\usepackage{amsmath, amsthm, amssymb, amsfonts, mathrsfs, amscd}
\usepackage{mathtools}
\usepackage[normalem]{ulem}
\usepackage{bm}
\usepackage{float} 
\usepackage{graphicx}
\usepackage{tikz-cd}
\usepackage{multirow, multicol}
\usepackage{longtable}
\usepackage{fouridx}
\usepackage{xcolor}
\usepackage[all, knot]{xy}
\usepackage{enumitem}
\usepackage[colorlinks=true, linkcolor=black, citecolor=black, urlcolor=blue, final]{hyperref}

\usepackage{cleveref}
\newcommand{\nc}{\newcommand}
\nc{\mfootnote}[1]{\footnote{#1}}
\nc{\todo}[1]{\textcolor{red}{[To do: #1]}}

\nc{\tred}[1]{\textcolor{black}{#1}}
\nc{\tblue}[1]{\textcolor{blue}{#1}}
\nc{\tgreen}[1]{\textcolor{green}{#1}}
\nc{\tpurple}[1]{\textcolor{purple}{#1}}
\nc{\btred}[1]{\textcolor{black}{\bfseries #1}}
\nc{\btblue}[1]{\textcolor{blue}{\bfseries #1}}
\nc{\btgreen}[1]{\textcolor{green}{\bfseries #1}}
\nc{\btpurple}[1]{\textcolor{purple}{\bfseries #1}}

\newtheorem{defn}{Definition}[section]
\newtheorem{definition}{Definition}[section]
\newtheorem{thm}{Theorem}[section]
\newtheorem{theorem}{Theorem}[section]
\newtheorem{cor}{Corollary}[section]
\newtheorem{prop}{Proposition}[section]

\newtheorem{remq}{Remark}[section]
\numberwithin{equation}{section}

\title[Non-abelian cohomology and extensions of Hom-algebras via the $\boldsymbol{\beta}$-Nijenhuis--Richardson bracket]{Non-abelian cohomology and extensions of Hom-algebras via the $\boldsymbol{\beta}$-Nijenhuis--Richardson bracket}
\author[N. Saadaoui]{Nejib Saadaoui}
\address{Laboratory of Mathematics and Applications LR17ES11, \\
	Higher Institute of Computer Science and Multimedia of Gabès, \\
	University of Gabès, Tunisia}
\email{najib.saadaoui@isimg.tn}

\begin{document}
		
		\maketitle
		
		\begin{abstract}
			This paper develops a cohomology theory for Hom-Leibniz algebras using the 
			$\beta$-Nijenhuis--Richardson bracket and applies it to classify non-abelian 
			extensions. We introduce left, and right versions of the bracket, 
			each defining a graded Lie algebra structure on the space of $\beta$-cochains. 
			The main result establishes that equivalence classes of split extensions of a 
			Hom-Leibniz algebra $L$ by $V$ are in bijection with the second cohomology 
			space $H^2(L,V)$, generalizing classical results from Lie and Leibniz algebra 
			theory. We characterize extensions explicitly through 2-cocycles $(\lambda_l, 
			\lambda_r, \theta)$ and provide complete classifications of low-dimensional cases.
		\end{abstract}
		
		

		\medskip

		\bigskip
		\setcounter{tocdepth}{2}  % يعرض حتى subsection (اختياري)
		\hypersetup{linkcolor=blue} % تأكد أن الروابط في الفهرس زرقاء
		\tableofcontents
		\addtocontents{toc}{\protect\setcounter{tocdepth}{2}} % ي
		\bigskip
		
		\section{Introduction}

		The cohomology theory of Lie algebras, introduced by Chevalley and Eilenberg in the 1940s, has proven to be a fundamental tool in understanding the structure and classification of Lie algebras. A central result in this theory is the correspondence between equivalence classes of abelian extensions and second cohomology groups. This relationship has been extensively studied and generalized to various algebraic structures, including Leibniz algebras, which were introduced by Loday as a non-commutative generalization of Lie algebras.
		
		In recent years, the notion of Hom-algebras has emerged as a natural framework for studying deformations of classical algebraic structures. A Hom-algebra is characterized by a twisting map that modifies the standard identities. Hom-Lie algebras, introduced by Hartwig, Larsson, and Silvestrov \cite{HartwigLarssonSilvestrov}, and Hom-Leibniz algebras, studied by Makhlouf, Silvestrov~\cite{MakhloufSilvestrov2008}, have attracted considerable attention due to their connections with quantum deformations and discrete modifications of differential calculus.
		
		While the cohomology theory for Hom-Lie algebras has been developed by Sheng~\cite{REpSheng}, 
		and for Hom-Leibniz algebras by Cheng and Su~\cite{ChenLeibniz}, 
		the existing works mainly concern abelian extensions and representations. 
		The theory of non-abelian extensions—where the kernel algebra $V$ has a non-trivial multiplication—
		remains largely unexplored within the Hom-algebra framework. 
		Such extensions arise naturally when studying algebras whose derived subalgebra is non-abelian, 
		and their systematic treatment requires new techniques that go beyond the classical representation-theoretic approach. 
		This gap motivates the present work, which aims to extend the cohomological framework to the non-abelian setting.
		
		The primary objective of this paper is to develop a comprehensive cohomology theory for non-abelian extensions of Hom-Leibniz algebras and related structures. Our approach is based on the construction of graded Lie algebra structures via generalized Nijenhuis--Richardson brackets adapted to the Hom-algebra framework. We introduce the $\beta$-Nijenhuis--Richardson bracket, which encodes both the algebraic structure and the twisting map, and show that it provides a unified framework for studying left, right, and symmetric Hom-Leibniz algebras.
		
		Our main results concern the construction of suitable brackets, the description of extensions by cohomology classes, and explicit examples.
		
		The paper is organized as follows. 
		Section~\ref{sec:prelim} reviews the basic definitions and establishes the notation used throughout the paper. 
		Section~\ref{sec:NR} develops the $\beta$-Nijenhuis--Richardson brackets for left, right, and symmetric Hom-Leibniz algebras, and constructs the corresponding cohomology theories. 
		Section~\ref{sec:deformation} applies this framework to study one-parameter formal deformations of Hom-algebras, showing how the $\beta$-Nijenhuis--Richardson bracket provides a unified cohomological approach. 
		Section~\ref{sec:extensions} introduces the notion of 2-cocycles with values in non-abelian algebras and establishes the main correspondence between extensions and cohomology. 
		Finally, Section~\ref{sec:application} applies these results to classify two-dimensional regular Hom-Leibniz algebras, providing an algorithmic approach to classification via extension theory.

% في نهاية المقدمة (Introduction
		%%%%%%%%%%%%%%%%%%%%%%%%%%
		%%%%%%%%%%%%%%%%%%%%%%%%%%%%%%%%%%%%%%%%%%
		\section{Preliminaries and Notation}\label{sec:prelim}
		This section recalls the basic definitions of Hom-Leibniz algebras and introduces the notation used throughout the paper.
		\begin{defn}\label{isomHLei}
			A \emph{Hom-algebra} is a triple $(L, [\cdot,\cdot], \alpha)$ consisting of 
		a vector space $L$, a bilinear map $[\cdot,\cdot]\colon L\times L\to L$, 
		and a linear self-map $\alpha\colon L\to L$ called the \emph{twisting map}.
		
			A Hom-algebra $(L, [\cdot,\cdot], \alpha)$ is said to be \emph{regular} if $\alpha$ is an automorphism, and \emph{multiplicative} if 
\begin{equation}
    \alpha([x,y]) = [\alpha(x), \alpha(y)] \quad \text{for all } x, y \in L.
\end{equation}
		\end{defn}
		\begin{defn}[\cite{CharHomLeib,MakhloufSilvestrov2008,SaadaouiBiHomLeibniz}]\label{def:Symmetric Hom L}
			A Hom-algebra $(L, [\cdot,\cdot], \alpha)$ is said to be a:
			\begin{enumerate}[label=(\roman*)]
				\item \emph{left Hom-Leibniz algebra} if
				\[
				[\alpha(x), [y,z]] = [[x,y], \alpha(z)] + [\alpha(y), [x,z]], \quad \forall x,y,z \in L;
				\]
				\item \emph{right Hom-Leibniz algebra} if
				\[
				[\alpha(x), [y,z]] = [[x,y], \alpha(z)] - [[x,z], \alpha(y)], \quad \forall x,y,z \in L;
				\]
				\item \emph{symmetric Hom-Leibniz algebra} if both identities above hold.
			\end{enumerate}
		\end{defn}
		Note that setting $\alpha = \mathrm{id}_L$ recovers the classical Leibniz algebra identities 
		introduced by Loday~\cite{Loday1993}.
		%%%%%%%%%%%%%%%%%%
        
       \smallskip

By interchanging the variables $x$ and $y$ in the left Hom-Leibniz identity, 
one easily obtains the following property.

\begin{prop}\label{prop:symmetric_leibniz_identity}
If $(L, [\cdot,\cdot], \alpha)$ is a left Hom-Leibniz algebra, then for all $a,b,c \in L$,
\begin{equation}
    [[b,c],\alpha(a)] = -[[c,b],\alpha(a)].
\end{equation}
\end{prop} 
 Similarly, we have the following.
\begin{prop}\label{prop:right_leibniz_identity}
If $(L, [\cdot,\cdot], \alpha)$ is a right Hom-Leibniz algebra, then for all $a,b,c \in L$,
\begin{equation}\label{prop:left_leibniz_identity}
    [\alpha(a), [b,c]] = -[\alpha(a), [b,c].
\end{equation}
\end{prop} 
        %%%%%%%%%%%%%%%%%%%%%%%%%%%%%%%%%%%%%
		\subsection{The symmetric algebra $S(V)$}\label{subsection: symmetric algebra}
		Let $V$ be a finite-dimensional vector space over a field $\mathbb{K}$. The symmetric algebra is defined as
		\[
		S(V) = \bigoplus_{k=0}^{\infty} S^k(V),
		\]
		where $S^k(V)$ denotes the $k$-th symmetric power of $V$, with $S^0(V) = \mathbb{K}$ and $S^1(V) = V$.  
		As a vector space, $S(V)$ is isomorphic to the polynomial ring $\mathbb{K}[x_1,\dots,x_n]$, where $n = \dim V$.  
		Once a basis $\{v_1,\dots,v_n\}$ of $V$ is fixed, this isomorphism identifies each vector $v_i$ with the indeterminate $x_i$.
		
		Any linear map $f: V \to W$ extends uniquely to a $\mathbb{K}$-algebra homomorphism $S(f): S(V) \to S(W)$, determined by
		\[
		S(f)(v_{i_1} \cdots v_{i_k}) = f(v_{i_1}) \cdots f(v_{i_k}).
		\]
		In particular, if $f$ is an isomorphism, then so is $S(f)$. For simplicity, we shall denote this extension again by $f$.
		
		We denote by
		\[
		C^{\bullet}(V) = \bigoplus_{k \geq 1} \operatorname{Hom}\big(S^{k}(V), V\big)
		\]
		the graded vector space whose component of degree $k-1$ consists of all linear maps $S^k(V) \to V$.

	\subsection{Group action}	
		
	For $X \in \mathrm{Hom}(W^{\otimes n}, W)$ and $f: V \to W$ an isomorphism, we define
	\[
	(f\cdot X)(a_1,\dots,a_n)=f^{-1}\!\big(X(f(a_1),\dots,f(a_n))\big),
	\]
	that is,
	\begin{equation}\label{eq:action}
		f\cdot X = f^{-1} \circ X \circ f^{\otimes n}.
	\end{equation}
	This defines a natural action of $\mathrm{GL}(W)$ on $\mathrm{Hom}(W^{\otimes n}, W)$ when $V = W$.
	
	Then, two Hom-algebras $(V,[\cdot,\cdot],\alpha)$ and $(W,[\cdot,\cdot]',\alpha')$ are isomorphic if there exists an  isomorphism $f: V \to W$  such that
	\[
	\alpha = f \cdot \alpha' \quad \text{and} \quad [\cdot,\cdot] = f \cdot [\cdot,\cdot]'.
	\]
		
	\subsection{Hom-cochains}
	Let $(L, [\cdot, \cdot]_L, \alpha_L)$ be a Hom-algebra and $(V, \alpha_V)$ a vector space equipped with a linear map $\alpha_V \colon V \to V$. For each integer $k \geq 1$, a linear map
	\[
	f \colon L^{\otimes k} \to V
	\]
	is called a \emph{$k$-Hom-cochain} if it satisfies the following \emph{twisting condition}:
	\[
	\alpha_V \circ f = f \circ \alpha_L^{\otimes k}.
	\]
	We denote by $C_\alpha^k(L,V)$ the space of all $k$-Hom-cochains and define
	\[
	C_\alpha(L,V) = \bigoplus_{k \geq 1} C_\alpha^k(L,V).
	\]
	This space is naturally a $\mathbb{Z}$-graded vector space, where the homogeneous component of degree $k-1$ is $C_\alpha^k(L,V)$.
	
	The relevant subspace of $C_\alpha^k(L,V)$ depends on the specific structure of $L$:
	\begin{itemize}
\item If $L$ is a left or right Hom-Leibniz algebra, we consider the full space $C_\alpha^k(L,V)$.
		\item If $L$ is a Hom-Lie algebra, we restrict to the subspace of \emph{alternating} cochains, denoted $C_{\alpha,\mathrm{alt}}^k(L,V)$.

  \item  If $L$ is a \emph{symmetric} Hom-Leibniz algebra,
  we fix a semi-alternating product $d \in C_{\alpha}^n(L,L)$, that is, a fixed $n$-linear map $d \colon L^{\otimes n} \to L$ satisfying the identity
\begin{equation}\label{eq:semiAlternating1}
\begin{split}
d\big(d(x_1, \dots, x_n),\, \alpha(x_{n+1}), \dots, \alpha(x_{2n-1})\big)
&= (-1)^{i-1}\,
d\big(\alpha(x_{n+1}), \dots, \alpha(x_{n+i-1}),\\
			&\qquad\,
d(x_1, \dots, x_n),\, 
\alpha(x_{n+i}), \dots, \alpha(x_{2n-1})\big),
\end{split}
\end{equation}
for all $x_1, \dots, x_{2n-1} \in L$ and for each $i \in \{2, \dots, n\}$.

We then define the space of \emph{semi-alternating cochains}, denoted by $C_{\alpha,\mathrm{alt}'}^k(L,L)$, as the subspace of $k$-Hom-cochains $f \in C_\alpha^k(L,L)$ satisfying the following two graded symmetry conditions with respect to the fixed $d$:

For all $x_1, \dots, x_{k+n-1} \in L$ and for all $i \in \{2, \dots, k\}$,
\begin{equation}\label{eq:semiAlternating2}
\begin{split}
f\big(d(x_1, \dots, x_n),\,\alpha (x_{n+1}), \dots,\alpha( x_{k+n-1})\big)
&= (-1)^{i-1}\,
f\big(\alpha (x_{n+1}), \dots,\alpha (x_{n+i-1}),\\&\qquad
d(x_1, \dots, x_n),\, 
\alpha (x_{n+i}), \dots,\alpha (x_{k+n-1})\big),
\end{split}
\end{equation}

and
\begin{equation}\label{eq:semiAlternating3}
\begin{split}
d\big(f(x_1, \dots, x_k),\, \alpha( x_{k+1}), \dots,\alpha (x_{k+n-1})\big)
&= (-1)^{i-1}\,
d\big(\alpha(x_{k+1}), \dots, \alpha(x_{k+i-1}),\\&\qquad
f(x_1, \dots, x_k),\, 
\alpha(x_{k+i}), \dots, \alpha(x_{k+n-1})\big).
\end{split}
\end{equation}
	
\end{itemize}

In the next section, we define the $\beta$-Nijenhuis--Richardson bracket, where $\beta$ denotes the twisting map on the direct sum $M = L \oplus V$. This bracket acts naturally on the $\beta$-cochain space $C_\beta(M) = C_\beta(M,M)$.
	
		%%%%%%%%%%%%%%%%%%%%%%%%%%%% 
		%%%%%%%%%%%%%%%ù
		%%%%%%%%%%%%%%%%%%%%%%%%%
		%%%%%%%%%%%%%%%%%
		%%%%%%%%%%
		%%%%%%%%%%%%%%%%%%%%%%%%%%%%%%
		%\section{$\beta$-Nijenhuis--Richardson Brackets}
		%	\label{sec:NR}
Throughout this paper, we adopt the following conventions:
\begin{itemize}
    \item The base field $\mathbb{K}$ is assumed to have characteristic zero.
    
    \item All Hom-algebras (Hom-Leibniz algebras, Hom-Lie algebras, etc.) are assumed to be \textbf{multiplicative}.
\end{itemize}
            		
		\section{$\beta$-Nijenhuis--Richardson Brackets}
		\label{sec:NR}
		The classical Nijenhuis--Richardson bracket \cite{NIJENHUIS,NIJENHUIS2} 
		has played a central role in the development of deformation and cohomology theories 
		for Lie algebras. In particular, Fialowski and Penkava~\cite{Alice} 
		used this bracket to describe the structure of cohomology and deformations 
		of (super) Lie algebras.
		
		Building on this approach, we extend these ideas to the Hom setting 
		by introducing $\beta$-twisted left and right Nijenhuis--Richardson brackets 
		adapted to Hom-Leibniz and Hom-Lie algebras. 
		These brackets endow the space of $\beta$-cochains 
		with a graded Lie algebra structure, leading to cohomology complexes 
		that naturally generalize the classical case.

		%%%%%%%%%%%%%%%%%%%%%%%%%%%%%%%%%%%%%%%%%%%%%%%%%%%%%%%%%%%%
		The following result provides a $\beta$-twisted analogue of the graded Lie algebra structure
		introduced in~\cite{NIJENHUIS2}.
		
		\begin{defn}\label{def:gradedLie}
			Let $[\cdot,\cdot] \colon C_{\beta}(M) \times C_{\beta}(M) \to C_{\beta}(M)$ be a bilinear map. 
			The pair $\left(C_{\beta}(M), [\cdot,\cdot]\right)$ is called a \emph{graded Lie algebra} if, 
			for all homogeneous elements $f \in C_{\beta}^{m}(M)$, $g \in C_{\beta}^{n}(M)$, and 
			$h \in C_{\beta}^{p}(M)$ (so that $\deg(f) = m-1$, $\deg(g) = n-1$, $\deg(h) = p-1$), 
			the following conditions hold:
			\begin{enumerate}[label=(\roman*)]
				\item \textbf{Degree compatibility:}
				\[
				\deg([f,g]) = \deg(f) + \deg(g).
				\]
				
				\item \textbf{Graded skew-symmetry:}
				\[
				[f,g] = -(-1)^{\deg(f)\deg(g)} [g,f].
				\]
				
				\item \textbf{Graded Jacobi identity:}
				\begin{equation}\label{eq:Graded Jacobi identity}
						(-1)^{\deg(f)\deg(h)} [f,[g,h]] 
					+ (-1)^{\deg(g)\deg(f)} [g,[h,f]] 
					+ (-1)^{\deg(h)\deg(g)} [h,[f,g]] = 0.
				\end{equation}
			\end{enumerate}
		\end{defn}
		
		\begin{prop}
			Let $\circ \colon C_{\beta}(M) \times C_{\beta}(M) \to C_{\beta}(M)$ be a bilinear map. 
			If the pair $\left(C_{\beta}(M), \circ\right)$ is a \emph{graded pre-Lie algebra}, that is, it satisfies
			\begin{equation}\label{eq:preLie}
				(f \circ g) \circ h - f \circ (g \circ h) 
				= (-1)^{(n-1)(p-1)} \left( (f \circ h) \circ g - f \circ (h \circ g) \right)
			\end{equation}
			for all $f \in C_{\beta}^{m}(M)$, $g \in C_{\beta}^{n}(M)$, $h \in C_{\beta}^{p}(M)$, 
			then the bracket defined by
			\begin{equation}\label{eq:bracketFromPreLie}
				[f, g] = f \circ g - (-1)^{(m-1)(n-1)} g \circ f
			\end{equation}
			endows $C_{\beta}(M)$ with a graded Lie algebra structure.
		\end{prop}
		
		\begin{proof}
			The argument follows the classical one in~\cite{NIJENHUIS}. 
			The bracket~\eqref{eq:bracketFromPreLie} automatically satisfies graded skew-symmetry, 
			and the graded Jacobi identity follows directly from the pre-Lie relation~\eqref{eq:preLie}.
		\end{proof}
		 Let \((C_\beta(M), [\cdot,\cdot])\) be  a graded-Lie algebra.  Let \( d \in C^2_\beta(M) \) be a $2$-Hom-cochain satisfying $[d,d]=0$. 
Define a sequence of linear maps $D_k \colon C_{\beta}^k(M) \to C_{\beta}^{k+1}(M)$ by
\begin{equation}\label{def:coboundary}
	D_k(f) = [d, f].
\end{equation}
Since \( d \in C_{\beta}^2(M) \) and \( f \in C_{\beta}^k(M) \), we have \( D_k(f) \in C_{\beta}^{k+1}(M) \); 
thus the composition \( D_{k+1} \circ D_k \) is well-defined. 
This leads to the following result.

\begin{prop}\label{prop:coboundary}
	The map $D_k \colon C_\beta^k(M) \to C_\beta^{k+1}(M)$ satisfies
	\[
	D_{k+1} \circ D_k = 0.
	\]
	Hence, $\left( \bigoplus_{k \geq 1} C_\beta^k(M), D_k \right)$ forms a cochain complex, 
	called the \emph{cohomology complex of the  Hom-algebra $(M, d, \beta)$}.
\end{prop}

\begin{proof}
	%Since \((C_\beta(M), [.,.])\) is a graded-Lie algebra, it follows from Proposition~\ref{eq:preLie} 
	%that the bracket \([\cdot,\cdot]_l\), as defined in~\eqref{eq:bracketNR}, endows \(C_\beta(M)\) 
	%with a graded Lie algebra structure.
	
	For \(f \in C_\beta^k(M)\), we compute:
	\[
	D_{k+1} \circ D_k(f) = [d, [d, f]].
	\]
	Applying the graded Jacobi identity in \((C_\beta(M), [\cdot,\cdot])\) and noting that 
	\(\deg(d) = 1\) and \(\deg(f) = k-1\), we obtain
	\[
	(-1)^{k-1}[d, [d, f]] - [d, [f, d]] + (-1)^{(k-1)} [f, [d, d]]=0.
	\]
	Since \([d, d] = 0\) and, by graded skew-symmetry, \([f, d] = -(-1)^{k-1} [d, f]\), this simplifies to
	\[
	2(-1)^{k-1}[d, [d, f]]=0.
	\]
    Since the base field $\mathbb{K}$ has characteristic zero, we conclude that $[d, [d, f]] = 0$, i.e.,
\[
D_{k+1} \circ D_k(f) = 0.
\]
\end{proof}
Let $\left( \bigoplus_{k \geq 1} C_\beta^k(M), D_k \right)$ be the corresponding cochain complex. 
A $k$-Hom-cochain $f \in C_{\beta}^{k}(M)$ is called a \emph{$k$-cocycle} if $D_k(f) = 0$, 
and a \emph{$k$-coboundary} if there exists $g \in C_{\beta}^{k-1}(M)$ such that $f = D_{k-1}(g)$. 

We denote by $Z_{\beta}^{k}(M) = \ker(D_k)$ the space of $k$-cocycles and by 
$B_{\beta}^{k}(M) = \operatorname{im}(D_{k-1})$ the space of $k$-coboundaries. 
Since $D_{k} \circ D_{k-1} = 0$, we have $B_{\beta}^{k}(M) \subset Z_{\beta}^{k}(M)$.

The $k$-th cohomology group of the Hom-algebra  $M$ is then defined as the quotient space:
\[
H_{\beta}^{k}(M) = Z_{\beta}^{k}(M) / B_{\beta}^{k}(M).
\]

		%%%%%%%%%%%%%%%%%%%%%%%%%%
		\subsection{The Left $\beta$-Nijenhuis--Richardson Bracket}	
		We now define a bilinear operation on $C_\beta(M)$ that endows it with the structure of a graded Lie algebra. 
		This operation is referred to as the \emph{left $\beta$-Nijenhuis--Richardson bracket}.
		
		Let $M$ be a vector space and $\beta \colon M \to M$ a linear map. 
		For integers $p, q \geq 2$, denote by $\mathrm{Sh}(p, q)$ the set of $(p,q)$-shuffles, 
		i.e., permutations $\sigma \in S_{p+q}$ satisfying
		\[
		\sigma(1) < \cdots < \sigma(p)
		\quad \text{and} \quad
		\sigma(p+1) < \cdots < \sigma(p+q).
		\]
		
		Define a map $l \colon \mathrm{Sh}(p, q) \to \{1, \dots, q\}$ by:
		\begin{itemize}
			\item $l(\sigma) = 1$ if $\sigma(p) < \sigma(p+1)$;
			\item $l(\sigma) = i$ if 
			\[
			\sigma(p+1) < \cdots < \sigma(p+i-1) < \sigma(p) < \sigma(p+i) < \cdots < \sigma(p+q).
			\]
		\end{itemize}
		
		For $f \in C_\beta^m(M)$ and $g \in C_\beta^n(M)$, define the composition $f \circ_l g \in C_\beta^{m+n-1}(M)$ by
		\begin{multline}\label{eq:rondNR}
			f \circ_l g(a_1, \dots, a_{m+n-1}) = 
			\sum_{\sigma \in \mathrm{Sh}(n, m-1)} (-1)^{l(\sigma)-1} \epsilon(\sigma) \\
			\times f\Big( 
			\beta^{n-1}(a_{\sigma(n+1)}), \dots, \beta^{n-1}(a_{\sigma(n+i-1)}),\,
			g(a_{\sigma(1)}, \dots, a_{\sigma(n)}),\,
			\beta^{n-1}(a_{\sigma(n+i)}), \dots, \beta^{n-1}(a_{\sigma(m+n-1)})
			\Big),
		\end{multline}
		where $i = l(\sigma)$ (that is, $i$ indicates the position of the term $g(a_{\sigma(1)}, \dots, a_{\sigma(n)})$ 
		inside the argument list of $f$), and $\epsilon(\sigma)$ denotes the signature of the shuffle $\sigma$.
		
		Using this composition, we define the bracket $[f,g]_l \in C_\beta^{m+n-1}(M)$ by
		\begin{equation}\label{eq:bracketNR}
			[f,g]_l = f \circ_l g - (-1)^{(m-1)(n-1)} g \circ_l f.
		\end{equation}
		The operation $[\cdot,\cdot]_l$ is called the \emph{left $\beta$-Nijenhuis--Richardson bracket} ($\beta$-NR bracket for short).
		\begin{prop}\label{prop:preLie}
			The pair $(C_\beta(M), \circ_l)$ defines a graded pre-Lie algebra.
		\end{prop}
		
		\begin{proof}
			We prove that the identity~\eqref{eq:preLie} holds for all homogeneous elements 
			$f \in C_\beta^k(M)$, $g \in C_\beta^l(M)$, and $h \in C_\beta^r(M)$.
			
			Let $\mathrm{Sh}(p,q,r)$ denote the set of $(p,q,r)$-unshuffles, i.e., permutations 
			$\sigma \in S_{p+q+r}$ satisfying
			\[
			\sigma(1) < \cdots < \sigma(p), \quad
			\sigma(p+1) < \cdots < \sigma(p+q), \quad
			\sigma(p+q+1) < \cdots < \sigma(p+q+r).
			\]
			
			We define two index maps 
			\[
			l' \colon \mathrm{Sh}(p,q,r) \to \{1, \dots, q\}
			\quad \text{and} \quad
			l'' \colon \mathrm{Sh}(p,q,r) \to \{1, \dots, r\}
			\]
			as follows:
			\begin{itemize}
				\item $l'(\sigma) = 1$ if $\sigma(p) < \sigma(p+1)$;
				\item $l'(\sigma) = i$ if 
				\[
				\sigma(p+1) < \cdots < \sigma(p+i-1) < \sigma(p) < \sigma(p+i);
				\]
				\item $l''(\sigma) = 1$ if $\sigma(p+q) < \sigma(p+q+1)$;
				\item $l''(\sigma) = j$ if 
				\[
				\sigma(p+q+1) < \cdots < \sigma(p+q+j-1) < \sigma(p+q) < \sigma(p+q+j).
				\]
			\end{itemize}
			
			Next, we compute the explicit form of the composition 
			$(f \circ_l g) \circ_l h \in C_\beta^{k+l+r-2}(M)$.
			For all $a_1, \dots, a_{k+l+r-2} \in M$, we have (where, for brevity, we write  $i=l'(\sigma)$ and $j=l''(\sigma)$):
			\begin{multline}\label{eq:fgcirc1}
				\big( (f \circ_l g) \circ_l h \big)(a_1, \dots, a_{k+l+r-2})\\
				= 
				\sum_{\substack{\sigma \in \mathrm{Sh}(r,l,k-2) \\ 1 \leq i< j \leq k-2}}
				(-1)^{i+j + l+1} \, \epsilon(\sigma) \,
				f\Big(
				\beta^{r+l-2}(a_{\sigma(r+l+1)}), \dots, 
				\beta^{r+l-2}(a_{\sigma(r+l+i-1)}), \\
				\beta^{r-1}\big( h(a_{\sigma(1)}, \dots, a_{\sigma(r)}) \big),\,
				\beta^{r+l-2}(a_{\sigma(r+l+i)}), \dots, 
				\beta^{r+l-2}(a_{\sigma(r+l+j-1)}), \\
			 g	\big(\beta^{r-1}(a_{\sigma(r+1)}), \dots,\beta^{r-1} (a_{\sigma(r+l)}) \big),\,
				\beta^{r+l-2}(a_{\sigma(r+l+j)}), \dots,
				\beta^{r+l-2}(a_{\sigma(r+l+k-2)})
				\Big)
			\end{multline}
			\begin{multline}\label{eq:Prelie2}
				+\sum_{\substack{\sigma \in \mathrm{Sh}(r,l-1,k-1) \\ 1 \leq i \leq l-1,\, 1\leq j \leq k-1}} 
				(-1)^{i+j} \epsilon(\sigma)
				f\Bigg( \beta^{r+l-2}(a_{\sigma(r+l)}), \dots, \beta^{r+l-2}(a_{\sigma(r+l+j-2)}), \\
				g\Big(\beta^{r-1}(a_{\sigma(r+1)}), \dots,
				\beta^{r-1}(a_{\sigma(r+i-1)}), h(a_{\sigma(1)}, \dots, a_{\sigma(r)}),\beta^{r-1}(a_{\sigma(r+i)}), \dots, \\
				\beta^{r-1}(a_{\sigma(r+l-1)}) \Big), \beta^{r+l-2}(a_{\sigma(r+l+j-1)}), \dots, 
				\beta^{r+l-2}(a_{\sigma(r+l+k-2)}) \Bigg)
			\end{multline}
			\begin{multline}\label{eq:fgcirc3}
				+ \sum_{\substack{\sigma \in \mathrm{Sh}(r,l,k-2) \\1\leq j<i \leq k-2}} 
				(-1)^{i+j+l} \epsilon(\sigma)  f\Big( \beta^{r+l-2}(a_{\sigma(r+l+1)}), \dots, \beta^{r+l-2}(a_{\sigma(r+l+j-2)}), \\
			 g	\big(\beta^{r-1}(a_{\sigma(r+1)}), \dots,\beta^{r-1}( a_{\sigma(r+l)}) \big), \beta^{r+l-2}(a_{\sigma(r+l+j-1)}), \dots, \beta^{r+l-2}(a_{\sigma(r+l+i-1)}), \\
				\beta^{r-1}\big( h(a_{\sigma(1)}, \dots, a_{\sigma(r)}) \big), \beta^{r+l-2}(a_{\sigma(r+l+i)}), \dots, \beta^{r+l-2}(a_{\sigma(r+l+k-2)}) \Big).
			\end{multline}
			
			We also have:
			\begin{multline*}
				f \circ_l (g \circ_l h)(a_1, \dots, a_{k+l+r-2}) = 
				\sum_{\substack{\sigma \in \mathrm{Sh}(r,l-1,k-1) \\1\leq i\leq k-1, 1\leq j\leq l-1}} (-1)^{i+j+r+1} \epsilon(\sigma) \\
				\times  f\Bigg(
				\beta^{l+r-2}(a_{\sigma(l+r)}), \dots, \beta^{l+r-2}(a_{\sigma(l+r+i-2)}), \,
				g\big( \beta^{r-1}(a_{\sigma(r+1)}), \dots, \beta^{r-1}(a_{\sigma(r+j-1)}), \\
				h(a_{\sigma(1)}, \dots, a_{\sigma(r)}), \beta^{r-1}(a_{\sigma(r+j)}), \dots, \beta^{r-1}(a_{\sigma(r+l-1)}) \big), \,
				\beta^{l+r-2}(a_{\sigma(l+r+i-1)}), \dots, \beta^{r+l-2}(a_{\sigma(l+r+k-2)}) \Bigg).
			\end{multline*}
			
			Then 
			\[
			f \circ_l (g \circ_l h) - (f \circ_l g) \circ_l h = \eqref{eq:fgcirc1} + \eqref{eq:fgcirc3}.
			\]
			
			By interchanging $g$ and $h$, we obtain the expression for the composition 
			$(f \circ_l h) \circ_l g \in C_\beta^{k+l+r-2}(M)$:
			\begin{multline}\label{eq:fhcirc}
				\big( (f \circ_l h) \circ_l g - f \circ_l (h \circ_l g) \big)(a_1, \dots, a_{k+l+r-2}) =
				\sum_{\substack{\sigma \in \mathrm{Sh}(l,r,k-2) \\ 1 \leq i< j\leq k-2}}
				(-1)^{i+j+r} \, \epsilon(\sigma) \\
				\times f\Big(
				\beta^{r+l-2}(a_{\sigma(r+l+1)}), \dots,
				\beta^{r+l-2}(a_{\sigma(r+l+j-1)}), \,
				 h\big(\beta^{l-1}(a_{\sigma(l+1)}), \dots, \beta^{l-1}(a_{\sigma(l+r)}) \big), \\
				\beta^{r+l-2}(a_{\sigma(r+l+j)}), \dots,
				\beta^{r+l-2}(a_{\sigma(r+l+i-1)}), \,
				\beta^{l-1}\big( g(a_{\sigma(1)}, \dots, a_{\sigma(l)}) \big), \\
				\beta^{r+l-2}(a_{\sigma(r+l+i)}), \dots,
				\beta^{r+l-2}(a_{\sigma(r+l+k-2)}) \Big)
			\end{multline}
			\begin{multline}\label{eq:fhcirc2}
				+
				\sum_{\substack{\sigma \in \mathrm{Sh}(l,r,k-2) \\ 1\leq  j<i\leq k-2}}
				(-1)^{i+j+r+1} \, \epsilon(\sigma) \,
				f\Big(
				\beta^{r+l-2}(a_{\sigma(r+l+1)}), \dots,
				\beta^{r+l-2}(a_{\sigma(r+l+j-1)}), \\
				\beta^{r-1}\big( g(a_{\sigma(1)}, \dots, a_{\sigma(l)}) \big), \,
				\beta^{r+l-2}(a_{\sigma(r+l+j)}), \dots,
				\beta^{r+l-2}(a_{\sigma(r+l+i-1)}), \\
				\beta^{l-1}\big( h(a_{\sigma(l+1)}, \dots, a_{\sigma(r+l)}) \big), \,
				\beta^{r+l-2}(a_{\sigma(r+l+i)}), \dots,
				\beta^{r+l-2}(a_{\sigma(r+l+k-2)}) \Big),
			\end{multline}
			where $\epsilon(\sigma)$ denotes the signature of the shuffle $\sigma$.
			
		Given $\sigma \in \mathrm{Sh}(l, r, k-2)$ , we define the permutation $\tau$ as:
		\[
		\tau =
		\begin{pmatrix}
			1 & \cdots & r & r+1 & \cdots & r+l \\
			l+1 & \cdots & r+l & 1 & \cdots & l
		\end{pmatrix}.
		\]
		Using $\tau$, we define $\sigma' \in \mathrm{Sh}(r, l, k-2)$ by:
		\[
		\sigma' = \sigma \circ \tau.
		\]
		This gives the relationship between the signs of the permutations:
		\[
		\epsilon(\sigma) = (-1)^{rl} \epsilon(\sigma').
		\]
		
	Next, since $g$ and $h$ are $\beta$-Hom-cochains, and by substituting $\sigma$ with $\sigma' \tau^{-1}$ in equation~\eqref{eq:fhcirc}, a direct comparison of the terms in~\eqref{eq:fhcirc}--\eqref{eq:fhcirc2} with those in~\eqref{eq:fgcirc1}--\eqref{eq:Prelie2} shows that:
	
		\[
		(-1)^{(r-1)(l-1)} \eqref{eq:fhcirc} =\eqref{eq:fgcirc1} ,
		\quad \text{and} \quad
		(-1)^{(r-1)(l-1)} \eqref{eq:fhcirc2} =\eqref{eq:fgcirc3} .
		\]
		
		This confirms that the pre-Lie identity~\eqref{eq:preLie} holds.
		\end{proof}

	\smallskip
	For $d \in C_\beta^2(M)$, one computes
	\begin{equation*}
		\frac{1}{2} [d,d]_l(a,b,c) 
		= d(d(a,b), \beta(c)) + d(\beta(b), d(a,c)) - d(\beta(a), d(b,c)).
	\end{equation*}
	Hence, $[d,d]_l = 0$ if and only if $(M, d, \beta)$ is a left Hom-Leibniz algebra.
	
	\smallskip
	Let $(M, d, \beta)$ be a left Hom-Leibniz algebra and 
let \( d \in C^2_\beta(M) \) be a \(2\)-\(\beta\)-cochain. 
Define the linear map \( D_k \colon C_{\beta}^k(M) \to C_{\beta}^{k+1}(M) \) by
\begin{equation}\label{def:coboundary}
    D_k(f) = [d, f]_l.
\end{equation}
Then, by Proposition~\ref{prop:coboundary}, the operator \( D_k \) defines a \emph{coboundary operator}, 
and the sequence \( \left( C_{\beta}^k(M), D_k \right)_{k \ge 1} \) forms a cohomology complex.

%%%%%%%%%%%%%%%%%%%%%%%%%%%%%%%%%%%%%%%%%%%%%%%%%%%%%%%%%%%%%%%%%%%%%%%%%%%%%%%%%%%%%%

We now provide an explicit formula for the coboundary operator \( D_k = [d, \cdot]_l \).
For any \( f \in C_{\beta}^k(M) \) and \( a_1, \dots, a_{k+1} \in M \), we have
\begin{multline}\label{eq:coboundary_explicit}
    D_k(f)(a_1, \dots, a_{k+1}) = [d, f]_l(a_1, \dots, a_{k+1}) \\
    =(-1)^{k+1} \sum_{s=1}^{k} (-1)^{s+1}\,
        d\!\big( \beta^{k-1}(a_s),\, f(a_1, \dots, \widehat{a_s}, \dots, a_{k+1}) \big) \\
    + d\!\big( f(a_1, \dots, a_k),\, \beta^{k-1}(a_{k+1}) \big) \\
    +(-1)^{k+1} \sum_{1 \leq s < t \leq s+1} (-1)^{s}\,
        f\!\big( \beta(a_1), \dots, \widehat{a_s}, \dots, d(a_s, a_t), \dots, \beta(a_{k+1}) \big),
\end{multline}
where the symbol $\widehat{a_s}$ indicates that the argument \( a_s \) is omitted from the list.

% This formula coincides with the coboundary operator introduced in~\cite{ChenLeibniz}.
This characterization extends naturally to the symmetric case, as follows.
\begin{prop}
A triple $(M, d, \beta)$ with $d \in C^2_\beta(M)$ is a symmetric 
Hom-Leibniz algebra  if and only if 
 $ d$ defines a $2$-semi alterning product satisfying the identity $[d, d]_l = 0$.
\end{prop}

Explicitly, the conditions that $ d$ is a 
$2$-semi alterning product
 and that $[d, d]_l = 0$ 
are respectively equivalent to the following two identities, valid for all $a, b, c \in M$:
\begin{align}
    d\big( \beta(a), d(b, c) \big) &= -\, d\big( d(b, c), \beta(a) \big), 
    \label{eq:beta-d-antisym}\\
    d\big( \beta(a), d(b, c) \big) &= d\big( d(a, b), \beta(c) \big) 
    + d\big( \beta(b), d(a, c) \big). \label{eq:sym-right}
\end{align}

The second cohomology group of a symmetric Hom-Leibniz algebra $(M, d, \beta)$
is obtained by restricting $D_k$ to the subspace $C^k_{\beta,\mathrm{alt}'}$.
Its explicit expression is given by the map 
$D_k \colon C_{\beta,\mathrm{alt'}}^k(M) \to C_{\beta,\mathrm{alt'}}^{k+1}(M)$ defined as
\begin{multline}\label{eq:coboundarySym}
    D_k(f)(a_1, \dots, a_{k+1}) = (-1)^{k+1}
    \sum_{s=1}^{k+1} (-1)^{s+1}\, d\Big( \beta^{k-1}(a_s),\, f(a_1, \dots, \widehat{a_s}, \dots, a_{k+1}) \Big) \\
    +(-1)^{k+1} \sum_{1 \leq s < t \leq k+1} (-1)^{s+t}\, 
        f\Big( d(a_s, a_t),\, \beta(a_1), \dots, \widehat{a_s}, \dots, \widehat{a_t}, \dots, \beta(a_{k+1}) \Big).
\end{multline}
%where $\widehat{a_s}$ indicates that $a_s$ is omitted.

Similarly, the second cohomology group of  Hom-Lie algebra $(M, d, \beta)$ 
is obtained by restricting $D_k$ to the subspace $C^k_{\alpha,\mathrm{alt}}$. 
Its explicit form is again given by~\eqref{eq:coboundarySym}. 
%which coincides with the expression introduced in~\cite{REpSheng}.

%for all $f \in \mathrm{Alt}_\beta^k(M)$ and $a_1, \dots, a_{k+1} \in M$.

			%%%%%%%%%%%%%%%%%%%%%%%%%%
		\subsection{The right $\beta$-Nijenhuis--Richardson Bracket}	
	%	\subsection{Cohomology of Right Hom-Leibniz Algebras}
		\label{subsec:rightCohomology}
		
		In this subsection, we develop a graded Lie algebra structure on the space of $\beta$-cochains associated with a right Hom-Leibniz algebra. The corresponding bracket is referred to as the \emph{right $\beta$-Nijenhuis--Richardson bracket}.

		Let $M$ be a vector space and let $\beta \colon M \to M$ be a linear map. For integers $p, q \geq 2$, denote by $\mathrm{Sh}(p, q)$ the set of $(p,q)$-shuffles. 

We define a map $r \colon \mathrm{Sh}(p, q) \to \{1, \dots, p\}$ by:
\begin{itemize}
	\item $r(\sigma) = 1$ if $\sigma(1)<\sigma(p+1)  $,
	\item $r(\sigma) = i$ if

	\[
	\sigma(p+1) < \cdots < \sigma(p+i) < \sigma(1) < \sigma(p+i+1) < \cdots < \sigma(p+q).
	\]

\end{itemize}

%%%%%%%%%%%%%%%%%%%%%%%%%%%%%%%%%%%%%%

For $f \in C_{\beta}^m(M)$ and $g \in C_{\beta}^n(M)$, we define the \emph{right $\beta$-insertion operator} $f \circ_r g \in C_{\beta}^{m+n-1}(M)$ by
\begin{multline}\label{eq:rondNRr}
	f \circ_r g(a_1, \dots, a_{m+n-1}) = 
	\sum_{\sigma \in \mathrm{Sh}(n, m-1)} (-1)^{r(\sigma)-1} \epsilon(\sigma) \cdot \\
	f\Big( \beta^{n-1}(a_{\sigma(n+1)}), \dots, \beta^{n-1}(a_{\sigma(n+i)}),\ 
	g(a_{\sigma(1)}, \dots, a_{\sigma(n)}),\ 
	\beta^{n-1}(a_{\sigma(n+i+1)}), \dots, \beta^{n-1}(a_{\sigma(n+m-1)}) \Big),
\end{multline}
where $i = r(\sigma)$ and $\epsilon(\sigma)$ denotes the signature of the shuffle $\sigma$.
%%%%%%%%%%%%%%%%%%%%%%%%%%%%%%%%%%%%%%%%%%%%%%%%%%%%%%%%%%%%%
		
%%%%%%%%%%%%%%%%%%%%%%%%%%%%%%%%%%%%%%	
		%%%%%%%%%%%%%%%%%%%%%%%%%%%%%%%%%%%%%%%%%%%%%%%%%%%%%%%%%%%%%
		
		The \emph{right $\beta$-Nijenhuis--Richardson bracket} $[f,g]_r \in C_{\beta}^{m+n-1}(M)$ is defined by
		\begin{equation}\label{eq:bracketNRr}
			[f,g]_r = f \circ_r g - (-1)^{(m-1)(n-1)} g \circ_r f.
		\end{equation}
        
			For $d \in C_\beta^2(M)$, we compute
		\begin{equation}\label{eq:ddr}
			\frac{1}{2} [d,d]_r(a,b,c) = d(d(a,b), \beta(c)) - d(d(a,c), \beta(b)) - d(\beta(a), d(b,c)).
		\end{equation}

		Therefore, the triple $(M, d, \beta)$ defines a right Hom-Leibniz algebra if and only if $[d,d]_r = 0$.
		
		Analogous to the left case, the coboundary operator for a right Hom-Leibniz algebra 
\((M, d, \beta)\) is defined by
\[
    D_k(f) = [d, f]_r.
\]
This operator gives rise to the second cohomology group of the right Hom-Leibniz algebra. 
Its explicit expression is given by
\begin{multline*}
    D_k(f)(a_1, \dots, a_{k+1}) = [d, f]_r(a_1, \dots, a_{k+1}) \\
    =(-1)^{k+1} d\!\big( \beta^{k-1}(a_1),\, f(a_2, \dots, a_{k+1}) \big)
    + (-1)^{k+1}\sum_{s=2}^{k+1} (-1)^{s}\,
        d\!\big( f(a_1, \dots, \widehat{a_s}, \dots, a_{k+1}),\, \beta^{k-1}(a_s) \big) \\
    +(-1)^{k} \sum_{1 \leq s < t \leq k+1} (-1)^{t}\,
        f\!\big( \beta(a_1), \ldots, d(a_s, a_t),\, 
        \beta(a_{s+1}), \ldots, \widehat{a_t},\ldots,\beta(a_{k+1}) \big),
\end{multline*}
where $\widehat{a_t}$ indicates that the argument \(a_s\) is omitted.

This expression coincides with the coboundary operator introduced in~\cite{ChenLeibniz}.		
		%\medskip
		%\noindent
		%\medskip
		%\noindent
		%\textbf{Notation.} In this paper, the bracket $[\cdot,\cdot]$ without subscript may refer to either the left or right $\beta$-Nijenhuis--Richardson bracket.

		%%%%%%%%%%%%%%%%%%%%

%%%%%%%%%%%%%%%%%%%%%%%%%%%%%%%%%		
%%%%%%%%%%%%%%%%%%%%%%%%%%%%%%%%%%%%%%%%%%%

	\section{Deformations of Hom-Algebras} \label{sec:deformation}
The cohomological theory of formal deformations for Hom-algebras was initiated by Makhlouf and Silvestrov in  \cite{MakhloufSilvestrov2010}, and further developed by Ammar, Ejbehi, and Makhlouf in \cite{AmmarEjbehiMakhlouf2011}. In this section, we extend this framework to Hom-Leibniz algebras using the $\beta$-Nijenhuis--Richardson bracket.

	%\subsection{One-Parameter Formal Deformations}
	Let $[\cdot,\cdot]$ denote the $\beta$-Nijenhuis--Richardson bracket (left or right)  and let $(M, d, \beta)$ be a Hom-algebra.
	
	\begin{definition}
		A \textbf{formal one-parameter deformation} of the Hom-algebra \(M\) is given by a family of bilinear maps
		\[
		d_t = d_0 + t\,d_1 + t^2\,d_2 + \cdots,
		\]
		where each $d_i \in C_{\beta}^2(M)$ is a $2$-$\beta$-cochain, and $d_0$ is the original bracket of \(M\). 
	\end{definition}
	%%%%%%%%%%%%%%%%%%%%%%%%%%%%
 Since $(M, d_t, \beta)$ is a Hom-algebra, the deformed bracket $d_t$ must satisfy
 \[
 [d_t, d_t] = 0,
 \]
 where $[\cdot,\cdot]$ denotes the $\beta$-Nijenhuis--Richardson bracket.  
 Expanding $d_t =\displaystyle \sum_{i \ge 0} t^i d_i$ and comparing the coefficients of $t^s$ gives the following infinite system:
 \begin{equation}\label{main_eq}
 	\sum_{i=0}^{s} [d_i, d_{s-i}] = 0, 
 	\qquad \text{for all } s = 0, 1, 2, \dots
 \end{equation}
 
 Analyzing each equation by degree of $t$, we obtain:
 \begin{itemize}
 	\item For $s = 0$: The condition $[d_0, d_0] = 0$ expresses the Hom-Jacobi identity for the original Hom-algebra $(M, d_0, \beta)$.
 	
 	\item For $s = 1$: We get $[d_0, d_1] = 0$, hence $d_1$ is a $2$-cocycle in the cochain complex $C^2_\beta(M)$.
 	
 	\item For $s = 2$: The relation $[d_0, d_2] + \tfrac{1}{2}[d_1, d_1] = 0$ shows that $[d_1, d_1]$ is a $3$-coboundary, i.e., $[d_1, d_1] \in B^3(M)$.
 	
 	\item For $s \ge 3$: The general condition can be written as
 	\[
 	2[d_0, d_s] + \Psi(d_1, \dots, d_{s-1}) = 0,
 	\]
 	where
 	\[
 	\Psi(d_1, \dots, d_{s-1}) = \sum_{i=1}^{s-1} [d_i, d_{s-i}].
 	\]
 	In particular, $\Psi(d_1, \dots, d_{s-1})$ defines a $3$-coboundary in $C_\beta^3(M)$.
 \end{itemize}

 %%%%%%%%%%%%%%%%%%%%%%%%%%%%%%%%%%%%%%%%%%%%%%%
Conversely,	assuming that \(\Psi(d_1, \dots, d_{s-1})\) is a \(3\)-coboundary (i.e., belongs to \(B^3(M)\)), there exists a \(d'_s \in C_{\beta}(M)\) such that \(\Psi(d_1, \dots, d_{s-1}) = [d_0, d'_s]\). The left-hand side then becomes \( [d_0, 2d_s + d'_s]=0\).
	This leads to the following result:
	
	\begin{prop} \label{prop:formal-defo}
Let $(d_i)_{i \in \mathbb{N}}$ be a sequence of $2$-$\beta$-cochains. Then $\displaystyle d_t = \sum_{i \geq 0} t^i d_i$ is a formal deformation of the Hom-algebra $M$ if and only if the following conditions hold:
\begin{enumerate}[label=(\roman*)]
 \item $[d_0,d_0] = 0$,
    \item $d_1 \in Z^2(M)$ is a $2$-cocycle on $M$.
    
    \item For all $s \geq 2$, the obstruction $\Psi(d_1, \dots, d_{s-1})$ is a $3$-coboundary (i.e., belongs to $B^3(M)$), and for any choice of $d'_s \in C^2_{\beta}(M)$ such that $\Psi(d_1, \dots, d_{s-1}) = [d_0, d'_s]$, the cochain $2d_s + d'_s$ is a $2$-cocycle on $M$.
\end{enumerate}
\end{prop}
	%%%%%%%%%%%%%%%%%%%%%%%%%%%%%%%%%%%%%%
	%%%%%%%%%%%%%%%%%%%%%%%%%%%%%%%%%%%
\subsection{Equivalent and Trivial Deformations}	
\begin{defn}
Let $M_t = (M, d_t, \beta)$ and $M'_t = (M, d'_t, \beta)$ be two formal deformations of $M$.
We say that the two deformations are \textbf{equivalent} if there exists a \textbf{formal isomorphism}
\[
\Phi_t \colon M[[t]] \to M[[t]],
\]
i.e., a $K[[t]]$-linear map that may be written in the form
\[
\Phi_t = \sum_{i \geq 0} \Phi_i t^i = \Phi_0 + \Phi_1 t + \Phi_2 t^2 + \cdots,
\]
where $\Phi_i \in C^1_\beta(M)$ and $\Phi_0 = \mathrm{Id}$, such that the following relation holds:
\begin{align}
d'_t &= \Phi_t \cdot d_t, \label{eq:morphism1} 
\end{align}
where "$\cdot$" denotes the action defined in \eqref{eq:action}.

A deformation $M_t$ is called \textbf{trivial} if it is equivalent to the original algebra $M$ (viewed as a constant deformation over $K[[t]]$).
\end{defn}
	By comparing the coefficients of \(t\) in equation \eqref{eq:morphism1}, we find: $d'_0=d$, and $d_1'=d_1 +\Phi_1.d_0$ Hence, 
	\[
	d'_1 = d_1 + [d_0, \Phi_1],
	\]
    where $\Phi_1 \in C^1_\beta(M)$ and $[\cdot,\cdot]$ denotes the $\beta$-Nijenhuis--Richardson bracket.This means that $d'_1$ and $d_1$ differ by a $2$-coboundary. 
	. This leads to the following proposition:
	
\begin{prop}
	There is a one-to-one correspondence between the elements of the second cohomology space $H^2_\beta(M)$ and the equivalence classes of infinitesimal deformations $d_t = d_0 + d_1 t$ (over $\mathbb{K}[[t]]/(t^2)$). Consequently, if $H^2_\beta(M) = 0$, then every infinitesimal deformation of $M$ is trivial.
\end{prop}
	\begin{remq}
\begin{enumerate}
    \item If $d_0$ is any element of \( C_{\beta}^2(M) \) and each $d_t$ belongs to \( C_{\beta}^2(M) \), 
    then we obtain a deformation of the Hom-Leibniz algebra \((M, d_0, \beta)\).
    
    \item If each $d_i \in C^2_{\beta,\mathrm{alt'}}(M)$, 
    then we obtain a deformation of the symmetric Hom-Leibniz algebra \((M, d_0, \beta)\).
    
    \item If each $d_i \in C^2_{\beta,\mathrm{alt}}(M)$, 
    then we obtain a deformation of the Hom-Lie algebra \((M, d_0, \beta)\).
\end{enumerate}
\end{remq}	
%%%%%%%%%%%%%%%%%%%%%%%%%%%		
			\section{Extensions and Representations of Hom-Algebras}\label{sec:extensions}
			The primary objective of this section is to extend certain notions and results regarding extensions of Lie algebras \cite{Alice}   to more general cases: Hom-Leibniz algebras and Hom-Lie algebras.
			%%%%%%%%%%%%%%%%%%%%%%
			\begin{defn}\label{defExt}
				Let $(L, \delta, \alpha)$, $(V, \mu, \alpha_V)$, and $(M, d, \alpha_M)$ be Hom-algebras, and let $i \colon V \to M$, $\pi \colon M \to L$ be morphisms of Hom-algebras. The following sequence of Hom-algebras is a short exact sequence if $\operatorname{Im}(i) = \ker(\pi)$, $i$ is injective, and $\pi$ is surjective:
				\[
				0 \longrightarrow (V, \mu, \alpha_V) \stackrel{i}{\longrightarrow} (M, d, \alpha_M) \stackrel{\pi}{\longrightarrow} (L, \delta, \alpha) \longrightarrow 0,
				\]
				where
\begin{equation}\label{extension9}
					\alpha_M \circ i = i \circ \alpha_V \quad \text{and} \quad \alpha \circ \pi = \pi \circ \alpha_M.
				\end{equation}
				In this case, we call $(M, d, \alpha_M)$ an extension of $(L, \delta, \alpha)$ by $(V, \mu, \alpha_V)$.
				
				Such an extension is called \textbf{split} if there exists a linear map $s \colon (L, \delta, \alpha) \to (M, d, \alpha_M)$ such that $\pi \circ s = \operatorname{id}_L$. Such a map $s$ is called a \textbf{section} of the extension.

			\end{defn}
			%%%%%%%%%%%%%%%%%%%%%
			\subsection{ Cohomology of Hom-Algebras}
			In Section~\ref{sec:NR}, we defined the cohomology of a Hom-algebra on itself (adjoint representation). In this section, we define a  cohomology of a Hom-algebra $L$ on another Hom-algebra $V$.
			%%%%%%%%%%%%%%%%%%%%%%%%%%
			
			Let $M = L \oplus V$ and $\beta = \alpha \oplus \alpha_V $.
For convenience, let
\[
C^{k,l}(L,V) = \operatorname{Hom}(L^k V^l, V),
\]
where $L^k V^l$ denotes the subspace of $M^{\otimes(k+l)}$ spanned by tensors with $k$ factors from $L$ and $l$ factors from $V$, in all possible orders given by shuffles $\sigma \in \operatorname{Sh}(k, l)$.

The symmetric counterpart is
\[
S^{k,l}(L,V) = \operatorname{Hom}(S^k(L) \otimes S^l(V), V).
\]

Let $[\cdot, \cdot]$ denote a $\beta$-Nijenhuis--Richardson bracket (left or right), and let $d \in C^2(M)$. We want to determine the conditions under which $(M, d, \beta)$ is a Hom-algebra and $V$ is an ideal of $M = L \oplus V$.

			If $V$ is an ideal of the Hom-algebra $L \oplus V$, then $[d, d] = 0$ and there exist
			\begin{align*}
				\delta &\in C^2(L), \\
				\lambda_l &\in C^{1,1}(LV, V), \\
				\lambda_r &\in C^{1,1}(VL, V), \\
				\theta &\in C^{2,0}(L, V), \\
				\mu &\in C^2(V)
			\end{align*}
			such that $d = \delta + \lambda_r + \lambda_l + \mu + \theta$. Hence,
			\begin{gather*}
				[\delta + \lambda_r + \lambda_l + \mu + \theta, \delta + \lambda_r + \lambda_l + \mu + \theta] = 0.
			\end{gather*}
			Therefore,
			\begin{multline*}
				\underbrace{[\delta, \delta]}_{=0} + \underbrace{2[\delta, \lambda_r + \lambda_l] + [\lambda_r + \lambda_l, \lambda_r + \lambda_l] + 2[\mu, \theta]}_{\in C^{2,1}} \\
				+ \underbrace{2[\mu, \lambda_r + \lambda_l]}_{\in C^{1,2}} + \underbrace{2[\delta + \lambda_r + \lambda_l, \theta]}_{\in C^{3,0}} + \underbrace{[\mu, \mu]}_{=0} = 0.
			\end{multline*}
			From this decomposition, we deduce the following conditions:
			\begin{align}
				& 2[\delta, \lambda_r + \lambda_l] + [\lambda_r + \lambda_l, \lambda_r + \lambda_l] + 2[\mu, \theta] = 0, \label{repOCT} \\
				& [\mu, \lambda_r + \lambda_l] = 0, \label{compatible} \\
				& [\delta + \lambda_r + \lambda_l, \theta] = 0: 
				\theta . \label{cocycle2}
			\end{align}
			
			Since $d$ is a $2$-$\beta$ cochain, we also have the multiplicativity condition:
			\begin{gather*}
				(\delta + \lambda_r + \lambda_l + \mu + \theta) \circ (\alpha + \alpha_V) = 
				(\alpha + \alpha_V) \circ (\delta + \lambda_r + \lambda_l + \mu + \theta).
			\end{gather*}
			This decomposes into the following component-wise conditions:
			\begin{align}
				\delta \circ \alpha &= \alpha \circ \delta, \label{multiplicative0} \\
				\lambda_r \circ (\alpha_V\otimes \alpha) &= \alpha_V \circ \lambda_r, \label{multiplicative1} \\
				\lambda_l \circ  (\alpha\otimes  \alpha_V) &= \alpha_V \circ \lambda_l, \label{multiplicative2} \\
				\mu \circ  \alpha_V &= \alpha_V \circ \mu, \label{multiplicative3} \\
				\theta \circ (\alpha\otimes \alpha) &= \alpha_V \circ \theta. \label{cocyDec}
			\end{align}
			This leads us to introduce the following definition and proposition.
			%By \eqref{multiplicative1}, \eqref{multiplicative2}, and \eqref{repOCT}, we have:
			%%%%%%%%%%%%%%%%%%%%%%%%
			\begin{defn} 
				A \textbf{$2$-cocycle }  of a Hom-algebra $(L, \delta, \alpha)$ with values in $(V,\mu,\alpha_V)$ is a triple $(\lambda_l,\lambda_r,\theta)$ that
				satisfies   \eqref{repOCT}, \eqref{compatible}, \eqref{cocycle2},  \eqref{multiplicative1}, \eqref{multiplicative2}, and \eqref{cocyDec}		
			\end{defn}
			%%%%%%%%%%%%%%%%%%%%%
			We summarize the main results in the following theorem.
\begin{theorem}\label{thm:cocycle}
Let $(L, \delta, \alpha)$ and $(V, \mu, \alpha_V)$ be Hom-algebras. Then $(L \oplus V, d, \alpha \oplus \alpha_V)$ is a Hom-algebra with $V$ as an ideal if and only if 
\[d = \delta + \lambda_r + \lambda_l + \mu + \theta,\]
where $(\lambda_l, \lambda_r, \theta)$ is a $2$-cocycle of $L$ with values in $V$.
\end{theorem}
As a direct consequence of Theorem \ref{thm:cocycle}, we obtain the following characterization of split extensions.
\begin{cor}\label{cor:split-extension}
Let $(L, \delta, \alpha)$ and $(V, \mu, \alpha_V)$ be Hom-algebras. A sequence
\[
E: 0 \longrightarrow (V, \mu, \alpha_V) \stackrel{i_0}{\longrightarrow} (L \oplus V, d, \alpha \oplus \alpha_V) \stackrel{\pi_0}{\longrightarrow} (L, \delta, \alpha) \longrightarrow 0,
\]
where $i_0$ and $\pi_0$ are the natural inclusion and projection maps, defines a split extension of $L$ by $V$ if and only if
$$d = \delta + \lambda_l + \lambda_r + \mu + \theta,$$ 
where $(\lambda_l, \lambda_r, \theta)$ is a $2$-cocycle of $L$ with values in $V$.
We call it the \textbf{standard split extension} of $L$ by $V$.
\end{cor}			%%%%%%%%%%%%%%%%%%%%%%%%
 The next proposition explicitly describes the 2-cocycle conditions for left Hom-Leibniz algebras.
			\begin{prop}\label{prop:2-cocycle-hom-leibniz}
				The triple $(\lambda_r, \lambda_l, \theta)$ is a 2-cocycle of a left  Hom-Leibniz algebra $L$ with values in $V$ if and only if the following conditions are satisfied: 
				\begin{align}
					\lambda_{r}(\alpha_V(v),\alpha(x))&=\alpha_V  \left( \lambda_{r}(v,x)\right)\label{prop:2-cocycle-hom-leibniz0} \\
					\lambda_{l}(\alpha(x),\alpha_V(v))&=\alpha_V \left(  \lambda_{l}(x,v)\right) \label{prop:2-cocycle-hom-leibniz1}\\
					\lambda_{l}\left(\alpha(x), \lambda_{l}(y,v)\right)	&=\lambda_{l}\left( \delta(x,y),\alpha_{V}(v)\right) +\lambda_{l}\left(\alpha(y), \lambda_{l}(x,v)\right)+\mu\left(\theta(x,y),\alpha_V(v) \right) ; \label{prop:2-cocycle-hom-leibniz2}\\
					\lambda_{l}\left(\alpha(x), \lambda_{r}(v,y)\right)&	=
					\lambda_{r}\left( \lambda_{l}(x,v), \alpha(y)         \right) 	
					+\lambda_{r}\left(\alpha_{V}(v), \delta(x,y)                \right)+\mu\left(\alpha_V(v),\theta(x,y) \right)\label{prop:2-cocycle-hom-leibniz3}\\
					\lambda_{r}\left(\alpha_{V}(v),\delta(x,y) 
					\right) 
					& = \lambda_{r}\left(\lambda_{r}(v,x),\alpha(y)  \right) 
					+\lambda_{l}\left(\alpha(x),\lambda_{r}(v,y)\right)-\mu\left(\alpha_V(v),\theta(x,y) \right).\label{prop:2-cocycle-hom-leibniz4} \\
					%%%%%%%%%%%%%%%%%%%%%%%%%%%
					\lambda_l\left(   \alpha(x),\mu(v,w)\right)&=\mu \left(  \alpha_V(v),\lambda_l(x,w) \right)+ \mu \left(  \lambda_l(x,v),\alpha_V(w) \right)\label{prop:2-cocycle-hom-leibniz5} \\
					%%%%%%%%%%%%%%%%%%
					\mu\left(\alpha_V(v),\lambda_l(x,w)\right)&=\mu\left(    \lambda_r(v,x),\alpha_V(w)\right) +\lambda_l\left(   \alpha(x),\mu(v,w)\right)  \label{prop:2-cocycle-hom-leibniz6}\\
					\mu\left(\alpha_V(v),\lambda_r(w,x)\right)&=\mu\left(  \alpha_V(w),  \lambda_r(v,x)\right) +\lambda_r\left(   \mu(v,w),\alpha(x)\right) \label{prop:2-cocycle-hom-leibniz7} \\ 
					%%%%%%%%%%%%%%%%%%%%%%%%%
					\theta(\alpha(x),\alpha(y))&=\alpha_{V}\left( \theta(x,y)\right);\label{prop:2-cocycle-hom-leibniz8}\\ 
					%%%%%%%%%%%%%%%%%%%%%%%%%%%
					\theta\left( \delta(x,y),\alpha(z)\right)
					&+\theta\left( \alpha(y),\delta(x,z)\right)
					-\theta\left( \alpha(x),\delta(y,z)\right)\nonumber\\
					%%%%%%%%%%%%%%%%%%%%%%%%%
					+\lambda_{r}\left( \theta(x,y),\alpha(z)\right) 
					&+\lambda_{l}\left(\alpha(y),\theta(x,z) \right)
					-\lambda_{l}\left(\alpha(x),\theta(y,z) \right) =0\label{prop:2-cocycle-hom-leibniz9}                           
				\end{align}
				for all $x,y,z \in L$ and $v,w \in V$.
			\end{prop}
			\begin{proof}Since $M$ is a left Hom-Leibniz algebra, we use in this proof the left $\beta$-Nijenhuis--Richardson bracket.
				Using \eqref{multiplicative1}, \eqref{multiplicative2} and \eqref{cocyDec}, we obtain \eqref{prop:2-cocycle-hom-leibniz0}, \eqref{prop:2-cocycle-hom-leibniz1} and \eqref{prop:2-cocycle-hom-leibniz8}, respectively.
				
				By \eqref{eq:rondNR}, for all $f,g\in C^2(M,M)$ we have
				\begin{align*}
					f\circ_lg(a,b,c)=f(g(a,b),\beta(c))+f(\beta(b),g(a,c))-f(\beta(a),g(b,c)).
				\end{align*}
				Note that the hom-cochains $\delta, \lambda_r, \theta, \mu$ vanish outside their domains of definition.
				
				Using \eqref{repOCT} with $(a,b,c)=(x,y,v)$, we obtain \eqref{prop:2-cocycle-hom-leibniz2}.
				
				Using \eqref{repOCT} with $(a,b,c)=(x,v,y)$, we obtain \eqref{prop:2-cocycle-hom-leibniz3}.
				
				Using \eqref{repOCT} with $(a,b,c)=(v,x,y)$, we obtain \eqref{prop:2-cocycle-hom-leibniz4}.
				
				Using \eqref{compatible} with $(a,b,c)=(x,v,w)$, we obtain \eqref{prop:2-cocycle-hom-leibniz5}.
				
				Using \eqref{compatible} with $(a,b,c)=(v,x,w)$, we obtain \eqref{prop:2-cocycle-hom-leibniz6}.
				
				Using \eqref{compatible} with $(a,b,c)=(v,w,x)$, we obtain \eqref{prop:2-cocycle-hom-leibniz7}.
				
				Using \eqref{cocycle2} with $(a,b,c)=(x,y,z)$, we obtain \eqref{prop:2-cocycle-hom-leibniz9}.
			\end{proof}

           Having established the cocycle conditions for left Hom-Leibniz algebras, we now extend this framework to the symmetric case.

			\begin{prop}
					
				The triple \((\lambda_l, \lambda_r, \theta)\) is a 2-cocycle of a symmetric Hom-Leibniz algebra \(L\) with values in \(V\) if and only if the following conditions are satisfied:
				\begin{gather}
	%\lambda_{l}\left(\alpha(x),\lambda_{l}(y,v) \right) = -\lambda_{l}\left(\alpha(x),\lambda_{r}(v,y) \right)\label{s0}\\
					%%%%%%%%%%%%%%%%%%%%%
					\lambda_{l}\left(\alpha(x),\lambda_{l}(y,v) \right) = -\lambda_{r}\left(\lambda_{l}(y,v),\alpha(x) \right)\label{s1}\\
					%%%%%%%%%%%%%%%%%%%%%%%%%%%%
					\lambda_{l}\left(\alpha(x),\lambda_{r}(v,y) \right) = -\lambda_{r}\left(\lambda_{r}(v,y),\alpha(x) \right)\label{s2}\\ 					
					%%%%%%%
					\lambda_{r}\left(\alpha_V(v),\delta(x,y) \right)
					+\mu(\alpha_V(v),\theta(x,y))	=-
					\lambda_{l}\left(\delta(x,y),\alpha_V(v) \right) -\mu(\theta(x,y),\alpha_V(v))\label{s3}\\
					%%%%%%%%%%%%%%%%%%%%%
					%%%%%%%
					\lambda_{l} \left(\alpha(x),\mu(u,v) \right) = -\lambda_{r}\left(\mu(u,v),\alpha(x) \right).\label{s6}  \\			%%%%%%%%%%%%%%%%%%%%%
					%=-\label{s7}\\% %%%%%%%
					\mu\left(\alpha_V(u),\lambda_{l}(x,v) \right) = -\mu\left(\lambda_{l}(x,v),\alpha_V(u) \right)\label{s5}\\
					%%%%%%%%%%%%%%%%%%% 
					\mu\left(\alpha_V(u),\lambda_{r}(v,x) \right) = -\mu\left(\lambda_{r}(v,x),\alpha_V(u) \right)\label{s4} 
					%%%%%%%%%%%%%%%%%%%%%%%%%%%%			
				\end{gather}
				and Equations \eqref{prop:2-cocycle-hom-leibniz1}, \eqref{prop:2-cocycle-hom-leibniz2}, \eqref{prop:2-cocycle-hom-leibniz5}, \eqref{prop:2-cocycle-hom-leibniz8}, and \eqref{prop:2-cocycle-hom-leibniz9}.
			\end{prop}        
			\begin{proof}
We apply equation \eqref{eq:beta-d-antisym} with $d = \delta + \lambda_r + \lambda_l + \theta + \mu$ and consider the following substitutions for $(a,b,c)$:
\begin{align*}
(a,b,c) = (x,y,v) &\quad \text{yields} \quad \eqref{s1}, \\
(a,b,c) = (x,v,y) &\quad \text{yields} \quad \eqref{s2}, \\
(a,b,c) = (v,x,y) &\quad \text{yields} \quad \eqref{s3}, \\
(a,b,c) = (x,u,v) &\quad \text{yields} \quad \eqref{s6}, \\
(a,b,c) = (u,x,v) &\quad \text{yields} \quad \eqref{s5}, \\
(a,b,c) = (u,v,x) &\quad \text{yields} \quad \eqref{s4}.
\end{align*}
%Substituting $(x,y,v)$ into \eqref{prop:left_leibniz_identity} gives \eqref{s0}.
Since every symmetric Hom-Leibniz algebra is in particular a left Hom-Leibniz algebra, the cocycle conditions \eqref{prop:2-cocycle-hom-leibniz0}--\eqref{prop:2-cocycle-hom-leibniz9} are satisfied.

Applying  \eqref{prop:left_leibniz_identity} with \((a,b,c)=(x,v,y)\) and using  \eqref{s1}, and \eqref{s3}, we deduce that equation \eqref{prop:2-cocycle-hom-leibniz3} leads to \eqref{prop:2-cocycle-hom-leibniz2}.

Similarly, Applying  \eqref{prop:left_leibniz_identity} with \((a,b,c)=(x,v,y)\) and 
using  \eqref{s3}, and \eqref{s2}, we deduce that equation \eqref{prop:2-cocycle-hom-leibniz4} leads to \eqref{prop:2-cocycle-hom-leibniz2}.

Applying \eqref{prop:left_leibniz_identity} with $(a,b,c) = (v,x,w)$, 
%and using \eqref{s5},
we find that \eqref{prop:2-cocycle-hom-leibniz6} leads to \eqref{prop:2-cocycle-hom-leibniz5}.

Finally, applying \eqref{prop:right_leibniz_identity} with $(a,b,c) = (v,x,w)$ and using \eqref{s4}, we find that \eqref{prop:2-cocycle-hom-leibniz7} leads to \eqref{prop:2-cocycle-hom-leibniz5}.

This completes the proof.
\end{proof}
 We now specialize the cocycle conditions to the case of Hom-Lie algebras, where antisymmetry imposes
			\begin{prop}\label{prop:2-cocycle-hom-lie}
The triple \( (\lambda_l, \lambda_r, \theta) \) is a \emph{2-cocycle} of a Hom-Lie algebra \( L \) with values in \( V \) if and only if it satisfies the following conditions:
\begin{align}
\lambda_r(v,x) &= -\lambda_l(x,v), \label{eq:lambda_relation} \\
\theta(x,y) &= -\theta(y,x), \label{eq:theta_relation}
\end{align}
together with the cocycle conditions \eqref{prop:2-cocycle-hom-leibniz1}, \eqref{prop:2-cocycle-hom-leibniz2}, \eqref{prop:2-cocycle-hom-leibniz5}, \eqref{prop:2-cocycle-hom-leibniz8}, and \eqref{prop:2-cocycle-hom-leibniz9}.
\end{prop}
\begin{proof}
Since the bracket of a Hom-Lie algebra is antisymmetric, we have $d(a,b) = -d(b,a)$ for all $a,b \in L \oplus V$. 
Setting $(a,b) = (x,v)$ with $x \in L$ and $v \in V$ yields \eqref{eq:lambda_relation}, 
and setting $(a,b) = (x,y)$ with $x,y \in L$ yields \eqref{eq:theta_relation}.

Moreover, every Hom-Lie algebra is a symmetric Hom-Leibniz algebra. 
Therefore, the general cocycle conditions for symmetric Hom-Leibniz algebras apply. 
Using the antisymmetry relations \eqref{eq:lambda_relation} and \eqref{eq:theta_relation}, 
one verifies that the remaining conditions reduce precisely to 
\eqref{prop:2-cocycle-hom-leibniz1}, \eqref{prop:2-cocycle-hom-leibniz2}, 
\eqref{prop:2-cocycle-hom-leibniz5}, \eqref{prop:2-cocycle-hom-leibniz8}, and \eqref{prop:2-cocycle-hom-leibniz9}.
\end{proof}
			\begin{remq}	
For the right Hom-Leibniz algebra, we obtain a similar result using the \emph{$\beta$-NR bracket}$ [\cdot, \cdot]_r$.
			\end{remq}
		%=================================
			\begin{prop}
				Let $f$ be a basis of a one-dimensional Hom-Leibniz algebra $L$, and let $v$ be a basis of a one-dimensional Hom-Leibniz algebra $V$. The set of 2-cocycles of $(L, \delta, \alpha)$ with values in $(V,\mu , \alpha_V)$ is described in the following table.
				\begin{table}[H]
					\centering
					\begin{tabular}{|c|c|c|c|c|c|}
\hline
\text{Hom-Leibniz algebra} & $\alpha_V(v) $ & $ \alpha(f)$ & $\lambda_l(f,v)$ & $ \lambda_r(v,f) $ & $\theta(f,f)$ \\
\hline
$L^0_1$ & $av$&  $bf$& -& -&-  \\
\hline
$L^1_{1}$ & $b^2v$&  $bf\ (b\neq1)$& -& -&$zv$ ($z\neq 0$)  \\
\hline
$L^2_{1}$ & $v$& $f$& $xv$& -& $zv$ ($z\neq 0$)  \\\hline
$L^3_{1}$ &$av$ & $f$&\(xv\) ($x\neq 0$) &- & -\\
\hline
$L^4_{1}$ &  $av$&$f$&$xv$ &$-xv$ & -\\
						\hline
				\end{tabular}\caption{Set of 2-cocycles for $L  $ on $V$}	\label{tab:cocycles}
				\end{table}
			\end{prop}
			\begin{proof}
Let $(f)$ be a basis of $L$ and $(v)$ be a basis of $V$. Set $\alpha_V(v)=av$, $\alpha(f)=bf$

Since $L$ and $V$ are one-dimensional regular left Hom-Leibniz algebras, we have $\delta=0$, $\mu=0$, and the twisting maps are invertible, so $a \neq 0$, $b \neq 0$.

Set  $\lambda_l(f,v)=xv$, $\lambda_r(v,f)=yv$, and $\theta(f,f)=zv$. Substituting into \eqref{prop:2-cocycle-hom-leibniz0}--\eqref{prop:2-cocycle-hom-leibniz9} yields the system:
\begin{equation*}
    (\mathcal{S})\quad
    \begin{cases}
        ya(b-1) = 0, \\
        xa(b-1) = 0, \\
        yb(y+x) = 0, \\
        z(b^2-a) = 0, \\
        zby = 0.
    \end{cases}
\end{equation*}
\textbf{Case 1: $b\neq 1$.} This case corresponds to the algebra $\mathbf{L^1_{1}}$.
Since $b\neq 1$ and $a \neq 0$, equations (1) and (2) give $x=y=0$. Then equation (5) is satisfied. Assuming $z\neq 0$ (for a non-trivial cocycle), equation (4) yields $b^2=a$.

\textbf{Case 2: $b=1$.} This case covers the algebras $\mathbf{L^0_{1}}$, $\mathbf{L^2_{1}}$, $\mathbf{L^3_{1}}$, and $\mathbf{L^4_{1}}$.
System $(\mathcal{S})$ reduces to $y(y+x)=0$, $z(1-a)=0$, and $zy=0$.
\item If $a=1$, the solutions are $(x,y,z) = (x,0,z)$ and $(x,-x,0)$.
    The solution $(x,0,z)$  corresponds to $\mathbf{L^2_{1}}$.
    The solution $(x,-x,0)$ corresponds to $\mathbf{L^4_{1}}$.
    \item If $a \neq 1$, then $z=0$. 
    The solution $(x,0,0)$ with $x \neq 0$ corresponds to $\mathbf{L^3_{1}}$,
    while the trivial cocycle $(0,0,0)$ corresponds to $\mathbf{L^0_{1}}$.
Thus, all cases in Table~\ref{tab:cocycles} are accounted for.
\end{proof}
%%%%%%%%%%%%%%%%%%%%%%%%%%%ù
\subsection{Equivalence of Extensions of Hom-Algebras}
			\label{subsec:equivalence}
In this subsection, we establish when two extensions are equivalent. 
The proofs make extensive use of the symmetric algebra framework 
(Subsection~\ref{subsection: symmetric algebra}): for elements $a, b$ 
in a vector space $M$, we write $ab$ for their product in $S(M)$, 
and bilinear maps are evaluated on these products.
			%%%%%%%%%%%%%%%%%%%%
			\begin{defn}\label{defEq}
				Two extensions 
				\[ 
				0 \longrightarrow (V, \mu, \alpha_{V}) \stackrel{i_1}{\longrightarrow} (M_1, d_1, \beta_1) \stackrel{\pi_1}{\longrightarrow} (L, \delta, \alpha) \longrightarrow 0, 
				\]
				\[ 
				0 \longrightarrow (V, \mu, \alpha_{V}) \stackrel{i_2}{\longrightarrow} (M_2, d_2, \beta_2) \stackrel{\pi_2}{\longrightarrow} (L, \delta, \alpha) \longrightarrow 0 
				\]
				are equivalent if there exists a Hom-algebra isomorphism  $\Phi \colon (M_1, d_1, \beta_1) \to (M_2, d_2, \beta_2)$  such that the following diagram commutes:
				\[
				\begin{tikzcd}
					0 \arrow[r] & (V, \mu, \alpha_{V}) \arrow[d, "\mathrm{id}_V"] \arrow[r, "i_1"] & (M_1, d_1, \beta_1) \arrow[d, "\Phi"] \arrow[r, "\pi_1"] & (L, \delta, \alpha) \arrow[d, "\mathrm{id}_L"] \arrow[r] & 0 \\
					0 \arrow[r] & (V, \mu, \alpha_{V}) \arrow[r, "i_2"] & (M_2, d_2, \beta_2) \arrow[r, "\pi_2"] & (L, \delta, \alpha) \arrow[r] & 0
				\end{tikzcd}
				\]
				that is,  $\Phi \circ i_1 = i_2 $ and $\pi_2 \circ \Phi = \pi_1$.
			\end{defn}
            \begin{remq}           
Equivalently, two extensions are equivalent if and only if there exists a linear isomorphism $\Phi \colon M_1 \to M_2$ such that $\Phi \cdot d_1 = d_2$, $\Phi \cdot \beta_1 = \beta_2$, $\Phi \circ i_1 = i_2$, and $\pi_2 \circ \Phi = \pi_1$.
\end{remq}
			%%%%%%%%%%%%%%%%%%%%%%%%
	\begin{theorem}\label{thm:surjective}
				Let 
				\[
				E: 0 \longrightarrow (V, \mu, \alpha_V) \stackrel{i}{\longrightarrow} (M, d', \alpha_M) \stackrel{\pi}{\longrightarrow} (L, \delta, \alpha) \longrightarrow 0,
				\]   
				be a split extension of $L$ by $V$. Then there exists a bilinear map $d = \delta + \lambda_l + \lambda_r + \theta + \mu$ such that:
				\begin{enumerate}[label=(\roman*)]
					\item $(\lambda_l, \lambda_r,\theta)$ is a  
					$2$-cocycle of $L$ with values in $V$,
					\item The short exact sequence
					\[
					E_0: 0 \longrightarrow (V, \mu, \alpha_V) \stackrel{i_0}{\longrightarrow} (L \oplus V, d, \alpha \oplus \alpha_V) \stackrel{\pi_0}{\longrightarrow} (L, \delta, \alpha) \longrightarrow 0,
					\]
					is an extension equivalent to $E$.
				\end{enumerate}
			\end{theorem}
			%%%%%%%%%%%
			\begin{proof}
				Let 
				\[
				0 \longrightarrow (V, \mu, \alpha_V) \stackrel{i}{\longrightarrow} (M, d', \alpha_M) \stackrel{\pi}{\longrightarrow} (L, \delta, \alpha) \longrightarrow 0,
				\]	
				be a split extension of $L$. Then there exists a section $s\colon L\to M$ such that $\pi \circ s=id_L$.
				Consider the isomorphism $\Phi : L \oplus V \rightarrow M$ defined by
                \[
	\Phi(x + v) = s(x) + i(v) \quad \text{for all } x \in L, v \in V.
	\]
 We have $ \Phi(\alpha(x) + \alpha_V(v) = s(\alpha(x)) + i(\alpha_V(v))=\alpha(s(x))+\alpha_V(i(v)) =(\alpha\oplus\alpha_V)(\Phi(x+v))$. 
				
				Define a bilinear map $d \colon L \oplus V \rightarrow L \oplus V$ by 
$d = \Phi^{-1} \cdot d'$.
				Then
\begin{equation}\label{ext2021}
	\Phi.[d,d]=[\Phi.d,\Phi.d]=[d',d']=0
\end{equation}
				Then $[d,d]=0$ and $(L \oplus V,d,\alpha\oplus\alpha_V)$
				is a Hom-algebra. Therefore, there exist a $2$-cocycle $(\lambda_l,\lambda_r,\theta)$ such that  $d=\delta+\lambda_r+\lambda_l+\theta+\mu$
				, and the short exact sequence
				\[
	E_0: 0 \longrightarrow (V, \mu, \alpha_V) \stackrel{i_0}{\longrightarrow} (L \oplus V, d, \alpha \oplus \alpha_V) \stackrel{\pi_0}{\longrightarrow} (L, \delta, \alpha) \longrightarrow 0,
				\]
				is a split extension of $L$ by $V$ equivalent to $E$. 
			\end{proof}
			As a consequence, every split extension of $(L, \delta, \alpha)$ by $(V, \mu, \alpha_V)$ is equivalent to the standard split extension $E_0$ defined above.
          %%%%%%%%%%%%%%%%%%%%%%%%%%%%%%%%ù
			\begin{prop}\label{prop:well def}
				Let two extensions \[E_0: 0\longrightarrow (V,\mu,\alpha_{V})\stackrel{i_{0}}{\longrightarrow} (L\oplus V,d,\alpha\oplus\alpha_{V})\stackrel{\pi_{0}}{\longrightarrow }(L,\delta,\alpha) \longrightarrow 0, \]
				\[E':0\longrightarrow (V,\mu,\alpha_{V})\stackrel{i_0}{\longrightarrow} (L\oplus V,d',\alpha\oplus\alpha_{V})\stackrel{\pi_0}{\longrightarrow }(L,\delta,\alpha) \longrightarrow 0\]			are equivalent if and only if there exists a $1$-Hom-cochain $h$ such that 
				\begin{align}
					%\delta' &= \delta, \\
					%\mu' &= \mu, \\       
					\lambda_l' &= \lambda_l + [\mu,h],\label{eq:equivalent cocycle1} \\
					\lambda_r' &= \lambda_r +[\mu,h] , \label{eq:equivalent cocycle2}\\
\theta' &=\theta+[\delta+\lambda_l+\lambda_r,h]+\frac{1}{2}[[\mu,h],h].\label{eq:equivalent cocycle3}
				\end{align}		
			\end{prop}
			\begin{proof}
				We assume that $ E_{0} $ and $ E' $ are equivalent.
				Then there exists an isomorphism of Hom-algebras
				\begin{align*}
					\Phi\colon&  \left( L\oplus V,d',\alpha\oplus\alpha_V\right)  \to \left( L\oplus V,d,\alpha\oplus\alpha_V\right)
				\end{align*}
such that the following diagram commutes:
		\begin{center}
		\begin{tikzcd}
			0 \arrow[r] & (V,\mu,\alpha_{V}) \arrow[r, "i_0"] \arrow[d, "{\text{id}_V}"] & (L \oplus V,d\text{'},\alpha \oplus\alpha_V) \arrow[r, "\pi_0"] \arrow[d, "\Phi"] & (L,\delta,\alpha) \arrow[r] \arrow[d, "{\text{id}_L}"] & 0 \\
						0 \arrow[r] & (V,\mu,\alpha_{V}) \arrow[r, "i_0"] & (L \oplus V,d,\alpha\oplus\alpha_V) \arrow[r, "\pi_0"] & (L,\delta,\alpha) \arrow[r] & 0
					\end{tikzcd}
				\end{center}
Then $\Phi(v)=v$. Let $\Phi(x)=s(x)+h(x) $ where $s(x)\in L$ and $h(x)\in V$. By $id_L\circ \pi_0=\pi_0\circ\Phi$ we  have $\Phi(x)=x+h(x)$, we denote $h(x+v)=h(x)$, then  $\Phi=1+h$. Since $\Phi$ is an isomorphism of Hom-algebra, we have  $d'=\Phi.d$ , i.e.,$
	d' = \Phi^{-1} \circ d \circ (\Phi \otimes \Phi).$ A direct computation then yields, for all $a, b \in L \oplus V$,
\begin{align*}
					d'(ab)&=\Phi^{-1}\circ d\circ \Phi(ab)\\&= \Phi^{-1}d\circ (\Phi(a)\Phi(b))\\
					&=\Phi^{-1}d(   (a+h(a))(b+h(b)))\\
					&=\Phi^{-1}d(ab+ah(b)+h(a)b+h(a)h(b))\\
	&=\Phi^{-1}(d(ab)+\lambda_l(ah(b))+\mu(ah(b))+\lambda_r(h(a)b))+\mu(h(a)b)+\mu(h(a)h(b))\\
    &=
	d(ab)-h(\delta(ab))+\lambda_l(ah(b)+\mu(ah(b))+\mu(h(a)b)
+\lambda_r(h(a)b))+\mu(h(a)h(b))\\
&=
d(ab)+[\delta,h](ab)+[\lambda_l,h](ab)+[\lambda_r,h](ab)+[\mu,h](ab)+\frac{1}{2}[[\mu,h],h](ab)
\end{align*}
	where we used that $\Phi^{-1}(x+ v) = x+ v - h(x)$ for $x \in L$, $v \in V$. 
    
    Consequently, the components of $d'$ are given by
				\begin{align*}
					%\delta' &= \delta, \\
					%\mu' &= \mu, \\                     
					\lambda_l' &= \lambda_l + [\mu,h],\\
					%\label{eq:equivalent cocycle1} \\
					\lambda_r' &= \lambda_r +[\mu,h] , \\
					%\label{eq:equivalent cocycle2}\\
\theta' &=\theta+[\delta+\lambda_l+\lambda_r,h]+\frac{1}{2}[[\mu,h],h].%\label{eq:equivalent cocycle3}
				\end{align*}
				%where $S(d\circ h)(a,b)=d\circ h(a,b)+d\circ h(b,a) $.
					Moreover, the condition $\Phi \circ (\alpha \oplus \alpha_V) = (\alpha \oplus \alpha_V) \circ \Phi$ implies
	\(
	h \circ \alpha = \alpha_V \circ h,\)
	so $h$ is a $1$-Hom-cochain.
		
Conversely, suppose there exists a $1$-Hom-cochain $h \colon L \to V$ satisfying the above relations.  
Define $\Phi = \mathrm{id}_{L} + h$.  
Then $\Phi$ is a linear isomorphism, satisfies $\Phi \circ(\alpha \oplus \alpha_V) = (\alpha \oplus \alpha_V)\circ  \Phi$, and fulfills $\Phi .d = d '$.  
Hence, $\Phi$ is an isomorphism of Hom-algebras, and the extensions $E_0$ and $E'$ are equivalent.
\end{proof}
	When $V$ is abelian (i.e., $\mu = 0$), the cohomology 
reduces to the classical representation cohomology. For the remainder of this section, 
our goal is to introduce and extend some well-known results originally established 
for abelian extensions to the case of non-abelian extensions of Hom-algebras, 
specifically in the context of Hom-Lie and Hom-Leibniz algebras.
			
Let $\beta=\alpha\oplus \alpha_V$.			Define a map 
\[
d^{1} \colon C_{\beta}^{1}(L,V) \longrightarrow 
\bigl( C_{\beta}^{1,1}(LV,V),\, C_{\beta}^{1,1}(VL,V),\, C^{2,0}_{\beta}(L,V) \bigr)
\]
by
\[
d^1 h = \bigl( [\mu,h],\, [\mu,h],\, [\delta,h] + \tfrac{1}{2}[[\mu,h],h] \bigr),
\]
where
 the first component $[\mu,h]\in C_{\beta}^{1,1}(LV,V)$ acts only on elements of the form $(x,v)\in L\times V$, the second component $[\mu,h]\in C_{\beta}^{1,1}(VL,V)$ acts only on elements of the form $(v,x)\in V\times L$. 
	\begin{prop}
The map $d^1$ is linear, and $\operatorname{Im} d^1 \subset Z^2(L,V)$, where $Z^2(L,V)$ denotes the space of all $2$-cocycles of $L$ with values in $V$.
\end{prop}
            \begin{proof}
			Let $h, h' \in C^1_\alpha(L,V)$ be $1$-Hom-cochains, and define the automorphisms
\[
\Phi(x+v) = x+ v + h(x), \quad \Phi'(x+v) = x+ v + h'(x) \quad (x \in L,\ v \in V).
\]
Then
\[
(\Phi \circ \Phi')(x+v) = x+ v + h(x) + h'(x),
\]
so the combined transformation corresponds to the cochain $h + h'$.

Let $d = \delta + \mu$ and set $d' = (\Phi \circ \Phi') \cdot d$.  
By the characterization of cocycles established in the previous proof, the components of $d'$ are
\begin{align*}
\lambda_l &= [\mu, h + h'], \\
\lambda_r &= [\mu, h + h'],  \\
\theta   &= [\delta, h + h'] + \tfrac{1}{2} [[\mu, h + h'], h + h'].
\end{align*}

On the other hand, let $d'' = \Phi' \cdot d$. Then $d''$ has components
\begin{align*}
\lambda_l' &= [\mu, h'], \\
\lambda_r' &= [\mu, h'], \\
\theta'    &= [\delta, h'] + \tfrac{1}{2} [[\mu, h'], h'],
\end{align*}
and applying $\Phi$ to $d''$ adds the contribution of $h$, yielding
\begin{align*}
\lambda_l &= [\mu, h] + [\mu, h'], \\
\lambda_r &= [\mu, h] + [\mu, h'], \\
\theta    &= [\delta, h] + [\delta, h'] + \tfrac{1}{2} [[\mu, h], h] + \tfrac{1}{2} [[\mu, h'], h'].
\end{align*}

Comparing both expressions, we obtain
\[
d^1(h + h') = d^1 h + d^1 h',
\]
which proves the linearity of $d^1$. We now verify that $d^1 h$ indeed lies in 
\[
\big( C^{1,1}_{\beta}(LV,V),\, C^{1,1}_{\beta}(VL,V),\, C^{2,0}_{\beta}(L,V) \big).
\]
%Let $\beta = \alpha \oplus \alpha_V$. Since $[\mu,h]$, $[\delta,h]$, and $\tfrac{1}{2}[[\mu,h],h]$ are all elements of $C^2_{\beta}(M)$, their restrictions to the respective subspaces $L  V$, $V  L$, and $L  L$ yield precisely the components of $d^1 h$. 
Hence, $d^1 h$ belongs to the stated product space.

				Finally, we show that $d^1 h$ is a $2$-cocycle.  
Let $d = \delta + \mu$. Since $(M, d, \beta)$ is a Hom-algebra and $\Phi = 1 + h$ is an automorphism of Hom-algebra of $M$, the bracket $d' = \Phi \cdot d$ also defines a Hom-algebra structure on $M$.  
By Theorem~\ref{thm:cocycle}, the corresponding triple $(\lambda_l, \lambda_r, \theta) = d^1 h$ is therefore a $2$-cocycle.
			\end{proof}
		Since $d^1$ is linear, it is straightforward to verify that the relation defined  on $Z^2(L,V)$ by
\[
f \sim g \iff f - g \in B^2(L,V), \quad \text{where } B^2(L,V) = \operatorname{im}(d^1),
\]
is an equivalence relation. The corresponding quotient space is then defined by
\[
H^2(L,V) = Z^2(L,V) / B^2(L,V).
\]
We now show that $H^2(L,V)$ classifies split extensions of $L$ by $V$, extending 
the classical result from the abelian to the non-abelian setting.

Denote by $\mathrm{Ext}(L,V)$ the set of equivalence classes of split extensions 
of $(L, \delta, \alpha)$ by $(V, \mu, \alpha_V)$. For an extension $\mathcal{E}$, 
we write $[\mathcal{E}]$ for its equivalence class, and for a $2$-cocycle 
$c = (\lambda_l, \lambda_r, \theta)$, we write $[c]$ for its class in $H^2(L,V)$.

\begin{theorem}\label{thm:ext_cohomology}
There exists a natural bijection between $H^2(L,V)$ and $\mathrm{Ext}(L,V)$, given by
\[
[c] \longmapsto [\mathcal{E}_c],
\]
where $\mathcal{E}_c$ denotes the standard split extension defined by the $2$-cocycle 
$c = (\lambda_l, \lambda_r, \theta)$ via the bracket 
$d_c = \delta + \lambda_l + \lambda_r + \theta + \mu$.
\end{theorem}
			\begin{proof}
				Define a map $\Phi \colon H^2(L,V) \to \mathrm{Ext}(L,V)$ by 
				$$[c] \mapsto [\mathcal{E}_c],$$
				where $\mathcal{E}_c$ denotes the extension corresponding to the 2-cocycle $c=(\lambda_l,\lambda_r,\theta)$.
				
				We show that $\Phi$ is a bijection by establishing the following three properties:
				
				\textbf{Well-definedness and injectivity:} By Proposition~\ref{prop:well def}, we deduce that $\Phi$ is well-defined, i.e., the map does not depend on the choice of representative in the cohomology class. By the the converse, we deduce that $\Phi$ is injective, i.e., distinct cohomology classes give rise to non-equivalent extensions.
				
				\textbf{Surjectivity:} By Theorem~\ref{thm:surjective}, we deduce that $\Phi$ is surjective, i.e., every equivalence class of extensions arises from some 2-cocycle.
				
                Therefore, $\Phi$ is a bijection.
			\end{proof}
%%%%%%%%%%%%%%%%%%%%%%%%%%%%%%%%%%%%%%%%%%%
			\section{An application}\label{sec:application}
Let $(M,d,\beta)$ be a Hom-algebra that is not simple. Suppose there exist complementary subspaces $L$ and $V$ such that $M = L \oplus V$, $\beta(L) \subseteq L$, $\beta(V) \subseteq V$, and $d(M,V)\subseteq V.$ We can decompose \( d \) as follows:
			\[
			d = \delta + \theta + \lambda_l + \lambda_r + \mu,
			\]
			where:
			\begin{itemize}
				\item \(\delta \colon L \times L \to L\) is the projection of \( d \) onto \( L \) (from \( L \times L \)),
				\item \(\theta \colon L \times L \to V\) is the projection of \( d \) onto \( V \) (from \( L \times L \)),
				\item \(\lambda_l \colon L \times V \to V\) is the projection of \( d \) onto \( V \) (from \( L \times V \)),
				\item \(\lambda_r \colon V \times L \to V\) is the projection of \( d \) onto \( V \) (from \( V \times L \)),
				\item \(\mu \colon V \times V \to V\) is the projection of \( d \) onto \( V \) (from \( V \times V \)).
			\end{itemize}
			It is easy to verify, using the previous section, that \((L, \delta, \alpha)\) and \((V, \mu, \alpha_V)\) are Hom-algebras, and that \((\lambda_l, \lambda_r, \theta)\) is a \(2\)-cocycle of \(L\) with values in \(V\).
			\begin{prop}
Every $1$-dimensional regular left Hom-Leibniz algebra is isomorphic to the algebra with basis $(v_1)$, trivial bracket $[v_1,v_1] = 0$, and twisting map $\alpha(v_1) = a v_1$ for some $a \in \mathbb{K}^\times$.
\end{prop}
\begin{proof}	
				The proof is a straightforward.	
			\end{proof}		
%Let $(f_1,\ldots,f_n)$ be a basis of $L$ and $(v_1,\ldots,v_m)$ be a basis of $V$. 
			\begin{thm}
				Every $2$-dimensional regular  left Hom-Leibniz algebra is isomorphic to one of the following pairwise non-isomorphic Hom-Leibniz  algebras.
\begin{table}[H]
\centering
\begin{tabular}{|c|c|c|c|c|}
\hline
\textbf{Name} &\textbf{Homomorphism} & \textbf{Multiplication}  & \textbf{Hom-Lie } & \textbf{Symmetric} \\
\hline
$L_{1,1}^{0}$ & $\alpha(e_1) = ae_1$, $\alpha(e_2) = be_2$ & -& Yes & Yes \\
						\hline
$L_{1,1}^{1}$ &$\alpha(e_1) = e_1$,$\alpha(e_2) =a e_2$  & $[e_1,e_2] = e_2$, $[e_2,e_1] =- e_2$ & Yes & Yes \\
\hline
						$L_{1,1}^{2}$ &$\alpha(e_1) =b e_1$ , $\alpha(e_2) =b^2 e_2$  & $[e_1,e_1]=e_2$&No & Yes \\
						\hline
$L_{1,1}^{3}$ &$\alpha(e_1) = e_1$,$\alpha(e_2) =a e_2$  & $[e_1,e_2] = e_2,$ & No & No \\
\hline
$L_{1,1}^{4}$ &$\alpha(e_1) = e_1$,$\alpha(e_2) = e_2$  & $[e_1,e_1] = e_2,\,[e_1,e_2] = ce_2,$ & No & No\\
\hline
$L_{1,0}^{5}$ &$\alpha(e_1) = e_1$, $\alpha(e_2) = e_1+e_2$  & - & Yes & Yes \\
\hline
$L_{1,0}^{6}$ &$\alpha(e_1) = e_1$, $\alpha(e_2) = e_1+e_2$  & $[e_1,e_2] =e_1$, $[e_2,e_1] =-e_1$ & Yes & Yes \\
\hline
$L_{1,0}^{7}$ &$\alpha(e_1) = e_1$, $\alpha(e_2) = e_1+e_2$  & $[e_2,e_2] =e_1$ & No & Yes \\						
						\hline
					\end{tabular}
					\caption{Classification of 2-dimensional regular   Hom-Leibniz algebras}
					\label{tab:hjj_classification}
				\end{table}                    
			\end{thm}
			\begin{proof}
				Let $(v,f)$ be  a Jordan-basis of a Hom-Leibniz algebra $(M,d,\beta)$. Then the matrix of $\beta$ has one of the following forms: $\begin{pmatrix}
					b&0\\0&a
				\end{pmatrix}$, $\begin{pmatrix}
					b&1\\0&b
				\end{pmatrix}$.
				We assume that $d(M,M)\neq M$. Then, $dim(d(M,M))=1$ since $M$ is non-abelian. Therefore, $V=d(M,M)$ is a one-dimensional Hom-algebra. Hence $\mu=0$. Let $L$ be a subalgebra of $M$ such that $L\oplus V=M$.
					\textbf{Case 1: $\beta = \begin{pmatrix} b & 0 \\ 0 & a \end{pmatrix}$.}  
			 Then $\beta(L)\subset L$ and $M$ is an extension of $L$ by $V$. Let \( f \) be a basis of \( L \) and \( v \) be a basis of \( V \).   Using Table~\ref{tab:cocycles}, we obtain the following:
   If $L = L^0_1$, the cocycle is $(0,0,0)$, yielding $L^0_{1,1}$.

	 If $L = L^1_1$, the cocycle is $(0,0,zv)$. Let us set $e_1 = f/z$, $e_2 = v/z$. In this basis, the bracket relations become $[e_1, e_1] = e_2$,  which corresponds to $L^2_{1,1}$.
     
 If $L = L^2_1$, the cocycle is $(x,0,zv)$.
 Let us set $e_1 = f/z$, $e_2 = v/z$.
In this basis, the bracket relations become $[e_1, e_1] = e_2$, $[e_1, e_2] = \frac{x}{z} e_2$, which corresponds to $L^4_{1,1}$.
     
 	If $L=L^3_{1}$, the cocycle is $(x,0,0)$. Let us set $e_1 = f/x$, $e_2 = v/x$.
In this basis, the bracket relations become $[e_1, e_2] = e_2$,  which corresponds to $L^3_{1,1}$.
    
 If $L = L^4_1$, the cocycle is $(x,-x,0)$. Let us set $e_1 = f/x$, $e_2 = v/x$.
In this basis, the bracket relations become $[e_1, e_2] = e_2$ and $[e_2, e_1] = -e_2$,  which corresponds to $L^1_{1,1}$.
 
\textbf{Case 2: $\beta = \begin{pmatrix} b & 1 \\ 0 & b \end{pmatrix}$.}  
The multiplicativity condition $\beta \circ d = d \circ (\beta \otimes \beta)$ forces $b = 1$. Direct verification of the Hom-Leibniz identity under this constraint yields the algebras $L^5_{1,0}$, $L^6_{1,0}$, and $L^7_{1,0}$.				
			\end{proof}	
					
			The above procedure provides an effective inductive algorithm for the classification of non-simple Hom-algebras of dimension \( s \). Having fixed a maximal ideal \( V \subset M \) with \( \dim V < s \), %the quotient \( \mathfrak{q} = \mathfrak{M}/V \) is necessarily simple. 
			By the induction hypothesis, both \( V \) and \(L\) are known; in particular, the case where \( V \) is simple poses no difficulty, as such algebras are already classified in lower dimensions. The algebra \( M \) is then reconstructed as an extension of \( L \) by \( V \). Explicitly, we determine all equivalence classes of Hom-Leibniz 2-cocycles \( [(\lambda_l,\lambda_r,\theta)] \in H^2(L, V) \), and for each class, we construct the corresponding bracket on \(M\cong V \oplus L \). Isomorphic extensions are identified using structural invariants, yielding all non-isomorphic non-simple Leibniz algebras of dimension \( s \).
			
			We note that when the maximal ideal coincides with the derived algebra, i.e. \( V = d(M,M)\), the associated structure of $L$ vanishes (\( \delta= 0 \)). If \( M \) is symmetric Hom-Leibniz algebra, then \( [M, M] \) is a Lie algebra. The systematic classification of such symmetric Leibniz algebras—and more generally, of algebras for which the derived subalgebra is a Lie algebra—will be the subject of a forthcoming paper.


\begin{thebibliography}{99}

\bibitem{AmmarEjbehiMakhlouf2011}
F. Ammar, Z. Ejbehi, and A. Makhlouf,
Cohomology and deformations of Hom-algebras,
\emph{J. Lie Theory} \textbf{21} (2011), 813--836.

\bibitem{CassasEXTENSIONS}
J. M. Casas and M. A. Insua,
On universal central extensions of Hom-Leibniz algebras,
\emph{J. Algebra Appl.} \textbf{13} (2014), 1450067.

\bibitem{CassasClssif}
J. M. Casas, M. Ladra, B. A. Omirov, and I. A. Karimjanov,
Classification of solvable Leibniz algebras with naturally graded filiform nilradical,
\emph{Linear Algebra Appl.} \textbf{43} (2013), 2973--3000.

\bibitem{ChenLeibniz}
Y. S. Cheng and Y. C. Su,
(Co)homology and universal central extension of Hom-Leibniz algebras,
\emph{Acta Math. Sin. (Engl. Ser.)} \textbf{27} (2011), 813--830.

\bibitem{Alice}
A. Fialowski and M. Penkava,
Extensions of (super) Lie algebras,
\emph{Commun. Contemp. Math.} \textbf{11} (2009), 709--737.

\bibitem{HartwigLarssonSilvestrov}
J. T. Hartwig, D. Larsson, and S. D. Silvestrov,
Deformations of Lie algebras using $\sigma$-derivations,
\emph{J. Algebra} \textbf{295} (2006), 314--361.

\bibitem{CharHomLeib}
A. Nourou Issa,
Some characterizations of Hom-Leibniz algebras,
\emph{Int. Electron. J. Algebra} \textbf{14} (2013), 1--13.

\bibitem{Loday1993}
J.-L. Loday,
Une version non commutative des algèbres de Lie: les algèbres de Leibniz,
\emph{Enseign. Math. (2)} \textbf{39} (1993), 269--293.

\bibitem{MakhloufSilvestrov2008}
A. Makhlouf and S. Silvestrov,
Hom-algebra structures,
\emph{J. Gen. Lie Theory Appl.} \textbf{2} (2008), 51--64.

\bibitem{MakhloufSilvestrov2010}
A. Makhlouf and S. Silvestrov,
Notes on 1-parameter formal deformations of Hom-associative and Hom-Lie algebras,
\emph{Forum Math.} \textbf{22} (2010), 715--739.

\bibitem{Geoffrey}
M. Geoffrey and Y. Gaywalee,
Leibniz algebras and Lie algebras,
\emph{SIGMA Symmetry Integrability Geom. Methods Appl.} \textbf{9} (2013), 063.

\bibitem{NIJENHUIS}
A. Nijenhuis and R. Richardson,
Cohomology and deformation in graded Lie algebras,
\emph{Bull. Amer. Math. Soc.} \textbf{72} (1966), 1--29.

\bibitem{NIJENHUIS2}
A. Nijenhuis and R. Richardson,
Deformations of Lie algebra structures,
\emph{J. Math. Mech.} \textbf{17} (1967), 89--105.

\bibitem{SaadaouiBiHomLeibniz}
N. Saadaoui,
(Bi)Hom--Leibniz algebras,
\emph{arXiv preprint}, arXiv:1810.07830 [math.RA], 2018.

\bibitem{REpSheng}
Y. Sheng,
Representations of Hom-Lie algebras,
\emph{Algebr. Represent. Theory} \textbf{15} (2012), 1081--1098.

\end{thebibliography}
	\end{document}